# Diagonalization of Operator functions by algebraic methods [1]

M. Stiefenhofer [2]

ABSTRACT. We give conditions for local diagonalization of an analytic operator family $L(\varepsilon)$ according to $L(\varepsilon) = \psi(\varepsilon) \cdot \Delta(\varepsilon) \cdot \phi^{-1}(\varepsilon)$ with diagonal operator polynomial $\Delta(\varepsilon)$ and analytic near identity bijections $\psi(\varepsilon)$ and $\phi(\varepsilon)$. The family $L(\varepsilon)$ is acting between real or complex Banach spaces $B$ and $\bar{B}$.

The basic assumption is given by stabilization of the Jordan chains at length $k$ in the sense that no root elements with finite rank above $k$ are allowed to exist. Jordan chains with infinite rank may appear. Decompositions of the linear spaces $B$ and $\bar{B}$ are constructed with corresponding subspaces assumed to be closed. These assumptions ensure finite pole order equal to $k$ of the generalized inverse $L^{-1}(\varepsilon)$ at $\varepsilon = 0$. The Smith form and smooth continuation of kernels and ranges of $L(\varepsilon)$ to appropriate limit spaces at $\varepsilon = 0$ arise immediately.

An algebraically oriented and self-contained approach is used, based on a recursion that allows for construction of power series solutions of $L(\varepsilon) \cdot b = 0$. The power series solutions are convergent, as soon as analyticity of $L(\varepsilon)$ and continuity of related projections are assumed.

## *1. Introduction*

Given an analytic family of linear operators $L(\varepsilon) = \sum_{i=0}^{\infty} \varepsilon^i \cdot L_i$, $L \in C^\omega(U, L[B, \bar{B}])$, $U$ open with $B, \bar{B}$ real or complex Banach spaces, $0 \in U \subset \mathbb{K} = \mathbb{R}, \mathbb{C}$ and $C^\omega(U, L[B, \bar{B}])$ denoting the set of analytic mappings from $U$ to the Banach space $L[B, \bar{B}]$ of bounded linear operators from $B$ to $\bar{B}$.

We give conditions for diagonalization of $L(\varepsilon)$ in the sense that there exist families of analytic, near identity transformations $\phi(\varepsilon) \in L[B, B]$ and $\psi(\varepsilon) \in L[\bar{B}, \bar{B}]$ satisfying

$$\psi^{-1}(\varepsilon) \cdot L(\varepsilon) \cdot \phi(\varepsilon) = \Delta(\varepsilon) \tag{1.1}$$

with $\phi(0) = I_B, \psi(0) = I_{\bar{B}}$ and a diagonal polynomial

$$\Delta(\varepsilon) = S_1 P_1 + \varepsilon \cdot S_2 P_2 + \cdots + \varepsilon^k \cdot S_{k+1} P_{k+1} \tag{1.2}$$

of degree $k$ where $S_i \in L[B, \bar{B}]$ are bounded linear mappings for each $i = 1, 2, \ldots, k+1$ such that the restricted mappings $S_{i|N_i^c} \in L[N_i^c, R_i]$ are bounded invertible mappings from the mutually disjoint closed subspaces $\{N_i^c\}_{i=1}^{k+1} \subset B$ to the mutually disjoint closed subspaces $\{R_i\}_{i=1}^{k+1} \subset \bar{B}$ and where the mappings $P_i \in L[B, N_i^c]$ for each $i = 1, 2, \ldots, k+1$ are bounded linear projections. Throughout the paper we assume $L(\varepsilon) \not\equiv 0$.

The subspaces $N_{i+1}^c, i = 1, \ldots, k$ are defined from root elements in the nullspace $N[L_0] := L_0^{-1}(\{0\}) = \{b \in B \mid L_0 b = 0\}$ of the operator $L_0 \in L[B, \bar{B}]$ using Jordan chains of length $i$ in the form $(b_{i-1} \cdots b_0) \in B^i$ for the family $\{L_i\}_{i \in \mathbb{N}_{-1}}$. The Jordan chain of length $i$ generated by the root

---

element $b_0 \neq 0$ must satisfy the condition $L_0 b_{i-1} + \cdots + L_{i-1} b_0 = 0$. If $b_0 \in N_{i+1}^c$, then the root element $b_0$ possesses no Jordan chain of length $i + 1$, i.e. $rank(b_0) = i$. For details concerning Jordan chains, see Definition 1 in section 3.

Then, for (1.1), (1.2) to be true, we have to assume that the maximal finite rank is limited by $k \geq 0$. Under these assumptions, direct sums of $B$ and $\bar{B}$ are constructed by

$$
\begin{array}{ccccccccc}
B & = & N_1^c & \oplus & N_2^c & \oplus & \cdots & \oplus & N_{k+1}^c & \oplus & N_{k+1} \\
 & & \downarrow \boxed{S_1} & & \downarrow \boxed{S_2} & & \cdots & & \downarrow \boxed{S_{k+1}} & & \\
\bar{B} & = & R_1 & \oplus & R_2 & \oplus & \cdots & \oplus & R_{k+1} & \oplus & R_{k+1}^c
\end{array}
\quad (1.3)
$$

with corresponding block diagonal operator matrix $\Delta(\varepsilon)$ given by

$$
\Delta(\varepsilon) = \begin{pmatrix} S_1 & \cdots & 0 & 0 \\ \vdots & \ddots & \vdots & \vdots \\ 0 & \cdots & \varepsilon^k \cdot S_{k+1} & 0 \\ 0 & \cdots & 0 & 0 \end{pmatrix} \begin{array}{c} R_1 \cdots R_{k+1} R_{k+1}^c \\ N_1^c \\ \vdots \\ N_{k+1}^c \\ N_{k+1} \end{array}. \quad (1.4)
$$

We have $N[L_0] = N_2^c \oplus \cdots \oplus N_{k+1}^c \oplus N_{k+1}$ with $N_{k+1}$ representing the nullspace of $\Delta(\varepsilon)$ for every $\varepsilon \in U \setminus \{0\}$, whereas $R_{k+1}^c$ denotes an arbitrary, but closed complement of $R_1 \oplus \cdots \oplus R_{k+1}$ in $\bar{B}$. The projections to $N_i^c$ and $R_i$ defined by (1.3) are assumed to be continuous.

Section 2 takes a first look at the literature. For motivation of the construction process of subspaces and associated operators in (1.3), we state an example in section 3 with stabilization of the Jordan chains at $k = 3$. After that, the general recursion concerning the direct sums in (1.3) and the coefficients $\{\phi_i\}_{i \in \mathbb{N}}$ and $\{\psi_i\}_{i \in \mathbb{N}}$ of the transformations $\phi(\varepsilon)$ and $\psi(\varepsilon)$ are given. The required results concerning Jordan chains are summarized in Lemma 1 of section 3. In sections 4 and 5, the results are successively proved and summarized in Theorems 1, 2 and 3. Finally, in section 6, we give some implications of the results, followed by a second look at the literature.

## 2. Comparison to the Literature

Diagonalization of analytic matrix functions $L(\varepsilon)$ at an isolated singularity can be found in [10, Theorem 1.2]. Concerning diagonalization of finite meromorphic operator functions with $L_0$ a Fredholm operator of arbitrary index, see [9, Theorem 11.6.4]. In the paper at hand, diagonalization in infinite dimensions is performed. We restrict to analytic operator functions $L(\varepsilon)$ (and operator power series) and mention that a meromorphic operator function $M(\varepsilon)$ turns into an analytic operator function by setting $L(\varepsilon) \coloneqq \varepsilon^p \cdot M(\varepsilon)$ and $p$ denoting the order of the pole of $M(\varepsilon)$. Then, our results may be extended from the analytic to the meromorphic case.

In [5], [6] and [11], the restriction to finite dimensions is also removed, whereas diagonalization of $L(\varepsilon)$ in the sense of (1.1) is not performed. Instead, in [6, Theorem 2.4] and [11, Theorem 3.9], the existence of a smooth generalized inverse $L^{-1}(\varepsilon)$ is characterized by stabilization of the Jordan chains at $k \geq 0$ ($k = \gamma$ = gradient) and closedness of the subspaces $N_{k+1} \subset B$ and $R_1 \oplus \cdots \oplus R_{k+1} \subset \bar{B}$ (in our terminology), combined with closedness of associated complements. Then, it is shown that these subspaces can analytically be continued to kernels and ranges of $L(\varepsilon), \varepsilon \neq 0$, allowing the construction of a generalized inverse $L^{-1}(\varepsilon)$ of $L(\varepsilon)$ with pole of order $k$.

The constructions in [5], [6] and [11] are working without the fine resolutions of $B$ and $\bar{B}$ in (1.3) given by $B = N_1^c \oplus \cdots \oplus N_{k+1}^c \oplus N_{k+1}$ and $\bar{B} = R_1 \oplus \cdots \oplus R_{k+1} \oplus R_{k+1}^c$. However, it is precisely this resolution that allows for diagonalization of $L(\varepsilon)$ with respect to powers of $\varepsilon$ by the diagonal operator polynomial $\Delta(\varepsilon) = S_1 P_1 + \varepsilon \cdot S_2 P_2 + \cdots + \varepsilon^k \cdot S_{k+1} P_{k+1}$. In addition, the Taylor expansions of the transformations $\phi(\varepsilon) = \sum_{i=0}^{\infty} \varepsilon^i \cdot \phi_i$ and $\psi(\varepsilon) = \sum_{i=0}^{\infty} \varepsilon^i \cdot \psi_i$ can easily be calculated by an explicit construction procedure. In [11], the investigation is strictly limited to generalized inverse operators. The diagonal operator polynomial $\Delta(\varepsilon)$ can be further factorized by $\Delta(\varepsilon) = S_P \cdot P(\varepsilon)$ allowing optionally factorization to the constant operator $S_P$ or to the Smith form $P(\varepsilon)$.

In summary, the direct sum decompositions proposed in this paper can be seen as supplementing the results in [5], [6] and [11]. Real and complex Banach spaces are treated without difference (no contour integrals) based on elementary analysis of formal power series, designed to solve the system of undetermined coefficients $\sum_{i+j=k} L_i \cdot b_j = 0$, $k \geq 0$. Assuming stabilization of the Jordan chains, this algebraically oriented approach is creating a Toeplitz matrix used to define the coefficients of the power series $\phi(\varepsilon) = \phi_0 + \varepsilon \cdot \phi_i + \varepsilon^2 \cdot \phi_2 + \cdots$ that is shown to satisfy a defining equation of power series type. Then, as soon as analyticity of $L(\varepsilon)$ and closedness of the subspaces in (1.3) are assumed, the defining equation turns into analyticity entailing analyticity of its solution $\phi(\varepsilon)$. No pre-transformation or extension of spaces is necessary in this approach.

If the leading operator $L_0$ is given by a Fredholm operator in Banach spaces, stabilization of Jordan chains and closedness of subspaces in (1.3) are automatically satisfied and the assumptions reduce to analyticity of $L(\varepsilon)$ in a deleted neighbourhood of the origin. In a finite dimensional setting, $L_0$ is a Fredholm operator. In an infinite-dimensional Hilbert space setting, stabilization of Jordan chains, closedness of the range spaces $R_i$, $i = 1, \ldots, k+1$ in (1.3) and analyticity of $L(\varepsilon)$ have to be assumed. Finally, when moving from Hilbert spaces to general Banach spaces, closedness of all subspaces in (1.3) needs to be added to the assumptions.

The coefficients of the power series $\phi(\varepsilon) = \phi_0 + \varepsilon \cdot \phi_1 + \varepsilon^2 \cdot \phi_2 + \cdots$ are used to first transform $L(\varepsilon)$ to upper triangularity, followed by the second transformation $\psi(\varepsilon) = \psi_0 + \varepsilon \cdot \psi_1 + \varepsilon^2 \cdot \psi_2 + \cdots$ to perform a standard process of back substitution that reduces the upper triangular system to diagonal form. This basic algebraic and self-contained approach seems not to be restricted to Banach or vector spaces, e.g. replacing the field $\mathbb{K}$ by a ring $\mathbb{A}$ might be interesting. Concerning the application to power series, the results are summarized in Theorems 1 and 2.

The recursion for construction of the subspaces in (1.3) was developed in [14, 295-300] for general Banach spaces. Independently, in [8] a comparable recursion was introduced for derivation of the direct sum $\bar{B} = R_1 \oplus \cdots \oplus R_{k+1}$. This investigation was continued and optimized in [12, Theorem 7.8.3] using Jordan chains in an explicit way. The approach in [8] and [12] is restricted to a family $L(\varepsilon)$ of Fredholm operators of index 0 with an isolated singularity at $\varepsilon = 0$. Then, kernels and ranges simplify to $N[L(\varepsilon)] \equiv \{0\}$ and $R[L(\varepsilon)] \equiv \bar{B}$ for $\varepsilon \neq 0$. In [12], real and complex Banach spaces are also treated in a parallel way.

Finally, in [3, section 3] a recursion is developed, which turns out to be comparable to the recursion from section 3.1 in case of a linear operator pencil $L(\varepsilon) = L_0 + \varepsilon \cdot L_1$ that stabilizes with $N_{k+1} = \{0\}$ and $R_{k+1}^c = \{0\}$. In [4, section 4] this recursion is further analyzed in case of non-termination of the iteration by use of the Lemma of Zorn.

### *3. Jordan chains, an Example and the Recursion*

Our main tool for construction of the subspaces in (1.3) is given by curves $b(\varepsilon) = \sum_{i=0}^{\infty} \varepsilon^i \cdot b_i$ in $B$ that are mapped to $\bar{B}$ by $L(\varepsilon) = \sum_{i=0}^{\infty} \varepsilon^i \cdot L_i$ according to

$$L(\varepsilon) \cdot b(\varepsilon) = \left( \varepsilon^{k-1} \cdots 1 \right) \cdot \begin{pmatrix} L_0 & \cdots & L_{k-1} \\ & \ddots & \vdots \\ & & L_0 \end{pmatrix} \cdot \begin{pmatrix} b_{k-1} \\ \vdots \\ b_0 \end{pmatrix} + \sum_{l=k}^{\infty} \varepsilon^l \cdot \sum_{i+j=l} L_i b_j \quad (3.1)$$

$$=: \left( \varepsilon^{k-1} \cdots 1 \right) \cdot \Delta^k \cdot (b_{k-1} \cdots b_0)^T + \sum_{l=k}^{\infty} \varepsilon^l \cdot \sum_{i+j=l} L_i b_j \ .$$

The following definitions are adapted versions from [12, Definition 7.2.1].

**Definition 1.** (i) For $b_0 \in B$, define the rank of $b_0$ according to

$$rank(b_0) := \begin{cases} 0, & \text{if } b_0 \notin N[L_0] \\ \sup \{ k \geq 1 \mid \exists (b_{k-1} \cdots b_0) \in B^k \text{ with } \Delta^k \cdot (b_{k-1} \cdots b_0)^T = 0 \}. \end{cases} \quad (3.2)$$

(ii) The family $L(\varepsilon)$ stabilizes at $k \geq 0$, if $rank(b_0) \in \{0, \cdots, k\} \cup \{\infty\}$ for all $b_0 \in B$ and there exists $b_0$ with $rank(b_0) = k$.

(iii) For $b_0 \in N[L_0]$ with $b_0 \neq 0$, the $k$-tuple $(b_{k-1} \cdots b_0) \in B^k, k \geq 1$ is called a Jordan chain of length $k$ with root element $b_0$, if it satisfies $\Delta^k \cdot (b_{k-1} \cdots b_0)^T = 0$.

In the next step, the recursion for diagonalization of $L(\varepsilon)$ is shown using an example with $k = 3$.

**Example 1.** Consider the polynomial matrix family

$$L(\varepsilon) = \begin{pmatrix} 1 & 0 & \varepsilon^3 \\ 0 & \varepsilon^2 & \varepsilon + \varepsilon^3 \\ \varepsilon^3 & 0 & \varepsilon^2 \end{pmatrix} = [e_1 \ 0 \ 0] + \varepsilon \cdot [0 \ 0 \ e_2] + \varepsilon^2 \cdot [0 \ e_2 \ e_3] + \varepsilon^3 \cdot [e_3 \ 0 \ e_1 + e_2] \quad (3.3)$$

with $L(\varepsilon) = L_0 + \varepsilon L_1 + \varepsilon^2 L_2 + \varepsilon^3 L_3$, $B = \bar{B} = \mathbb{R}^3$ and $e_1, e_2, e_3 \in \mathbb{R}^3$ denoting standard basis vectors in $\mathbb{R}^3$. Then, the diagonal operator polynomial $\Delta(\varepsilon) = S_1 P_1 + \varepsilon \cdot S_2 P_2 + \cdots + \varepsilon^k \cdot S_{k+1} P_{k+1}$ from (1.2) with corresponding subspaces in (1.3) is constructed along the following lines.

**Stage 1.** Set $S_1 := L_0 = [e_1 \ 0 \ 0]$ with kernel $N_1 := N[S_1] = sp(e_2, e_3)$ and range $R_1 := R[S_1] = sp(e_1)$. Here $sp(\cdot)$ means span of vectors in brackets. Hence, first decompositions of $B$ and $\bar{B}$ are obtained by $B = sp(e_1) \oplus sp(e_2, e_3) =: N_1^c \oplus N_1$ with projection $P_1 = [e_1 \ 0 \ 0]$ to $N_1^c$ and $\bar{B} = sp(e_1) \oplus sp(e_2, e_3) =: R_1 \oplus R_1^c$ with projection $\mathcal{P}_1 = [e_1 \ 0 \ 0]$ to $R_1$. The operator $\left[S_{1|N_0}\right]^{-1} \in L[R_1, N_1^c]$ is represented by the matrix $S_1^+ = [e_1 \ 0 \ 0] \in L[\mathbb{R}^3]$. The first summand of $\Delta(\varepsilon)$ reads $S_1 P_1 = [e_1 \ 0 \ 0]$. Finally, define $E_1 = (E_{1,1}) = M_1 = (M_{1,1}) := I_B = [e_1 \ e_2 \ e_3]$ and note that all Jordan chains of length 1 may be written as $b_0 = M_{1,1} \cdot n_1, n_1 \in N_1, n_1 \neq 0$.

**Stage 2.** Set $\bar{S}_2 := L_1 \cdot M_1 = [0 \ 0 \ e_2]$ and $S_2 := (I_{\bar{B}} - \mathcal{P}_1) \cdot \bar{S}_2 = [0 \ 0 \ e_2]$ with $S_{2|N_1} = S_2 \cdot (I_B - P_1) = [0 \ 0 \ e_2]$. Now, using $S_2$, the subspaces $N_1$ and $R_1^c$ are further split according to $N_2 := N[S_{2|N_1}] = sp(e_2)$, $R_2 := R[S_{2|N_1}] = sp(e_2)$, implying $N_1 = sp(e_3) \oplus sp(e_2) =: N_2^c \oplus N_2$ and $R_1^c = sp(e_2) \oplus sp(e_3) =: R_2 \oplus R_2^c$ with associated projection $P_2 = [0 \ 0 \ e_3]$ to $N_2^c$ and $\mathcal{P}_2 = [0 \ e_2 \ 0]$ to $R_2$. The operator $\left[S_{2|N_1}\right]^{-1} \in L[R_2, N_2^c]$ is represented by the matrix $S_2^+ = [0 \ 0 \ e_2] \in L[\mathbb{R}^3]$. The second summand of $\Delta(\varepsilon)$ is given by $\varepsilon \cdot S_2 P_2 = \varepsilon \cdot [0 \ 0 \ e_2]$. Now, define $E_2$ by setting

$$E_2 = \begin{pmatrix} E_{1,2} \\ E_{2,2} \end{pmatrix} := \begin{pmatrix} I_B & S_1^{-1} \mathcal{P}_1 \cdot \bar{S}_2 \\ 0 & I_B \end{pmatrix}^{-1} \cdot \begin{pmatrix} 0 \\ I_B \end{pmatrix}$$

or equivalently as the solution of the equation

$$\begin{pmatrix} I_B & S_1^{-1}\mathcal{P}_1 \cdot \bar{S}_2 \\ 0 & I_B \end{pmatrix} \cdot \begin{pmatrix} E_{1,2} \\ E_{2,2} \end{pmatrix} = \begin{pmatrix} 0 \\ I_B \end{pmatrix}$$

and hence show that in case of (3.3)

$$E_2 = \begin{pmatrix} 0 \\ I_B \end{pmatrix} \in \mathbb{R}^{6\times 3} \quad \text{as well as} \quad M_2 = \begin{pmatrix} M_{1,2} \\ M_{2,2} \end{pmatrix} := \begin{pmatrix} I_B & 0 \\ 0 & M_1 \end{pmatrix} \cdot E_2 = \begin{pmatrix} 0 \\ I_B \end{pmatrix} \in \mathbb{R}^{6\times 3}.$$

All Jordan chains of length 2 are given by

$$\begin{pmatrix} b_1 \\ b_0 \end{pmatrix} = \begin{pmatrix} M_{1,1} & M_{1,2} \\ 0 & M_{2,2} \end{pmatrix} \cdot \begin{pmatrix} n_1 \\ n_2 \end{pmatrix} = \begin{pmatrix} n_1 \\ n_2 \end{pmatrix}, \quad n_2 \neq 0.$$

**Stage 3.** Set $\bar{S}_3 := [L_1\ L_2] \cdot M_2 = [0\ e_2\ e_3]$ and $S_3 := (I_{\bar{B}} - \mathcal{P}_1 - \mathcal{P}_2) \cdot \bar{S}_3 = [0\ 0\ e_3]$ with $S_{3|N_2} = S_3 \cdot (I_B - P_1 - P_2) = [0\ 0\ 0]$. Using $S_3$, the subspaces $N_2$ and $R_2^c$ are next split according to $N_3 := N[S_{3|N_2}] = sp(e_2)$, $R_3 := R[S_{3|N_2}] = \{0\}$, yielding $N_2 = \{0\} \oplus sp(e_2) =: N_3^c \oplus N_3$ and $R_2^c = \{0\} \oplus sp(e_3) =: R_3 \oplus R_3^c$ with projection $P_3 = [0\ 0\ 0]$ to $N_3^c$ and $\mathcal{P}_3 = [0\ 0\ 0]$ to $R_3$. The operator $[S_{3|N_2}]^{-1} \in L[R_3, N_3^c]$ can be represented by the matrix $S_3^+ = [0\ 0\ 0] \in L[\mathbb{R}^3]$. The third summand of $\Delta(\varepsilon)$ reads $\varepsilon^2 \cdot S_3 P_3 = \varepsilon^2 \cdot [0\ 0\ 0]$. Define $E_3$ by setting

$$E_3 = \begin{pmatrix} E_{1,3} \\ E_{2,3} \\ E_{3,3} \end{pmatrix} := \begin{pmatrix} I_B & S_1^{-1}\mathcal{P}_1 \cdot \bar{S}_2 & S_1^{-1}\mathcal{P}_1 \cdot \bar{S}_3 \\ 0 & I_B & S_2^{-1}\mathcal{P}_2 \cdot \bar{S}_3 \\ 0 & 0 & I_B \end{pmatrix}^{-1} \cdot \begin{pmatrix} 0 \\ 0 \\ I_B \end{pmatrix}$$

or equivalently as the solution of the equation

$$\begin{pmatrix} I_B & S_1^{-1}\mathcal{P}_1 \cdot \bar{S}_2 & S_1^{-1}\mathcal{P}_1 \cdot \bar{S}_3 \\ 0 & I_B & S_2^{-1}\mathcal{P}_2 \cdot \bar{S}_3 \\ 0 & 0 & I_B \end{pmatrix} \cdot \begin{pmatrix} E_{1,3} \\ E_{2,3} \\ E_{3,3} \end{pmatrix} = \begin{pmatrix} 0 \\ 0 \\ I_B \end{pmatrix}$$

and hence calculate in case of (3.3)

$$E_3 = \begin{pmatrix} 0 & 0 & 0 \\ 0 & -e_3 & 0 \\ e_1 & e_2 & e_3 \end{pmatrix} \in \mathbb{R}^{9\times 3}.$$

Now define

$$M_3 = \begin{pmatrix} M_{1,3} \\ M_{2,3} \\ M_{3,3} \end{pmatrix} := \begin{pmatrix} I_B & 0 & 0 \\ 0 & M_{1,1} & M_{1,2} \\ 0 & 0 & M_{2,2} \end{pmatrix} \cdot \begin{pmatrix} E_{1,3} \\ E_{2,3} \\ E_{3,3} \end{pmatrix} = \begin{pmatrix} 0 & 0 & 0 \\ 0 & -e_3 & 0 \\ e_1 & e_2 & e_3 \end{pmatrix} \in \mathbb{R}^{9\times 3}$$

with all Jordan chains of length 3 given by

$$\begin{pmatrix} b_2 \\ b_1 \\ b_0 \end{pmatrix} = \begin{pmatrix} M_{1,1} & M_{1,2} & M_{1,3} \\ 0 & M_{2,2} & M_{2,3} \\ 0 & 0 & M_{3,3} \end{pmatrix} \cdot \begin{pmatrix} n_1 \\ n_2 \\ n_3 \end{pmatrix} = \begin{pmatrix} n_1 \\ n_2 - n_{2,3} \cdot e_3 \\ n_3 \end{pmatrix}, \quad n_3 \neq 0.$$

Here $n_{2,3}$ denotes the second component of $n_3 \in N_3 = sp(e_2)$.

**Stage 4.** Define $\bar{S}_4 := [L_1\ L_2\ L_3] \cdot M_3 = [e_3, -e_3, e_1 + e_2]$ and $S_4 := (I_{\bar{B}} - \mathcal{P}_1 - \mathcal{P}_2 - \mathcal{P}_3) \cdot \bar{S}_4 = [e_3, -e_3, 0]$, implying $S_{4|N_3} = S_4 \cdot (I_B - P_1 - P_2 - P_3) = [0, -e_3, 0]$, $N_4 := N[S_{4|N_3}] = \{0\}$, $R_4 :=$ $R[S_{4|N_3}] = sp(e_3)$, as well as $N_3 = sp(e_2) \oplus \{0\} =: N_4^c \oplus N_4$ and $R_3^c = sp(e_3) \oplus \{0\} =: R_4 \oplus R_4^c$ with projection $P_4 = [0\ e_2\ 0]$ to $N_4^c$ and $\mathcal{P}_4 = [0\ 0\ e_3]$ to $R_4$. The fourth summand of $\Delta(\varepsilon)$ reads $\varepsilon^3 \cdot S_4 P_4 = \varepsilon^3 \cdot [0, -e_3, 0]$.

Summing up the previous steps, decomposition (1.3) of $B = \mathbb{R}^3$ and $\bar{B} = \mathbb{R}^3$ arises as follows.

$$
\begin{array}{ccccccccc}
& \overbrace{sp(e_1)}^{} & & \overbrace{sp(e_3)}^{} & & \overbrace{\{0\}}^{} & & \overbrace{sp(e_2)}^{} & & \overbrace{\{0\}}^{} \\
B = & N_1^c & \oplus & N_2^c & \oplus & N_3^c & \oplus & N_4^c & \oplus & N_4 \\
& \downarrow \boxed{S_1 = [e_1\ 0\ 0]} & & \downarrow \boxed{S_2 = [0\ 0\ e_2]} & & \downarrow \boxed{S_3 = [0\ 0\ e_3]} & & \downarrow \boxed{S_4 = [e_3\ -e_3\ 0]} & & \\
\bar{B} = & R_1 & \oplus & R_2 & \oplus & R_3 & \oplus & R_4 & \oplus & R_4^c \\
& \underbrace{sp(e_1)}_{} & & \underbrace{sp(e_2)}_{} & & \underbrace{\{0\}}_{} & & \underbrace{sp(e_3)}_{} & & \underbrace{\{0\}}_{}
\end{array} \quad (3.4)
$$

In particular, by $N_4 = \{0\}$, no Jordan chains of length 4 exist, and hence, the family $L(\varepsilon)$ stabilizes at $k = 3$. We can now calculate the operator polynomial

$$\Delta(\varepsilon) = S_1 P_1 + \varepsilon \cdot S_2 P_2 + \varepsilon^2 \cdot S_3 P_3 + \varepsilon^3 \cdot S_4 P_4 = [\,e_1, -\varepsilon^3 \cdot e_3, \varepsilon \cdot e_2\,] = \begin{pmatrix} 1 & 0 & 0 \\ 0 & 0 & \varepsilon \\ 0 & -\varepsilon^3 & 0 \end{pmatrix}.$$

Representing $\Delta(\varepsilon)$ by the decompositions in (3.4) implies diagonalization according to

$$\Delta(\varepsilon) = \begin{pmatrix} S_1 P_1 & 0 & 0 & 0 \\ 0 & \varepsilon \cdot S_2 P_2 & 0 & 0 \\ 0 & 0 & \varepsilon^2 \cdot S_3 P_3 & 0 \\ 0 & 0 & 0 & \varepsilon^3 \cdot S_4 P_4 \end{pmatrix} = \begin{pmatrix} [e_1\ 0\ 0] & 0 & 0 & 0 \\ 0 & \varepsilon \cdot [0\ 0\ e_2] & 0 & 0 \\ 0 & 0 & 0 & 0 \\ 0 & 0 & 0 & \varepsilon^3 \cdot [0\ -e_3\ 0] \end{pmatrix}.$$

Finally, we define

$$E_4 = \begin{pmatrix} E_{1,4} \\ E_{2,4} \\ E_{3,4} \\ E_{4,4} \end{pmatrix} := \begin{pmatrix} I_B & S_1^{-1} \mathcal{P}_1 \cdot \bar{S}_2 & S_1^{-1} \mathcal{P}_1 \cdot \bar{S}_3 & S_1^{-1} \mathcal{P}_1 \cdot \bar{S}_4 \\ 0 & I_B & S_2^{-1} \mathcal{P}_2 \cdot \bar{S}_3 & S_2^{-1} \mathcal{P}_2 \cdot \bar{S}_4 \\ 0 & 0 & I_B & S_3^{-1} \mathcal{P}_3 \cdot \bar{S}_4 \\ 0 & 0 & 0 & I_B \end{pmatrix}^{-1} \cdot \begin{pmatrix} 0 \\ 0 \\ 0 \\ I_B \end{pmatrix}$$

and hence calculate in case of (3.3)

$$E_4 = \begin{pmatrix} 0 & 0 & -e_1 \\ 0 & 0 & -e_3 \\ 0 & 0 & 0 \\ e_1 & e_2 & e_3 \end{pmatrix} \in \mathbb{R}^{12 \times 3},$$

as well as

$$M_4 = \begin{pmatrix} M_{1,4} \\ M_{2,4} \\ M_{3,4} \\ M_{4,4} \end{pmatrix} := \begin{pmatrix} I_B & 0 & 0 & 0 \\ 0 & M_{1,1} & M_{1,2} & M_{1,3} \\ 0 & 0 & M_{2,2} & M_{2,3} \\ 0 & 0 & 0 & M_{3,3} \end{pmatrix} \cdot \begin{pmatrix} E_{1,4} \\ E_{2,4} \\ E_{3,4} \\ E_{4,4} \end{pmatrix} = \begin{pmatrix} 0 & 0 & -e_1 \\ 0 & 0 & -e_3 \\ 0 & -e_3 & 0 \\ e_1 & e_2 & e_3 \end{pmatrix} \in \mathbb{R}^{12 \times 3}.$$

We can now calculate the constant coefficients $\phi_0 = M_{4,4} = [e_1\ e_2\ e_3] = I_B$ and $\psi_0 = S_1 \cdot S_1^{-1} \mathcal{P}_1 + \cdots + S_4 \cdot S_4^{-1} \mathcal{P}_4 = \mathcal{P}_1 + \cdots + \mathcal{P}_4 = I_{\bar{B}}$ of the transformations $\phi(\varepsilon)$ and $\psi(\varepsilon)$ in (1.1). Although the decomposition is complete, the process can be continued indefinitely, above all for calculation of the coefficients $\phi_i, \psi_i, i \in \mathbb{N}$.

**Stage 5.** Define $\bar{S}_5 := [L_1 \; L_2 \; L_3 \; L_4] \cdot M_4 = [0, -(e_1 + e_2), e_3]$ and calculate $S_5 := (I_{\bar{B}} - \mathcal{P}_1 - \mathcal{P}_2 - \mathcal{P}_3 - \mathcal{P}_4) \cdot \bar{S}_5 = [0 \; 0 \; 0]$. Along the same lines as in Stages 1 to 4, we obtain subsequently

$$E_5 = \begin{pmatrix} 0 & e_1 & 0 \\ 0 & e_3 & 0 \\ 0 & 0 & 0 \\ 0 & 0 & -e_2 \\ e_1 & e_2 & e_3 \end{pmatrix} \in \mathbb{R}^{15 \times 3} \quad \text{and} \quad M_5 = \begin{pmatrix} 0 & e_1 & 0 \\ 0 & e_3 & -e_1 \\ 0 & 0 & 0 \\ 0 & -e_3 & -e_2 \\ e_1 & e_2 & e_3 \end{pmatrix} \in \mathbb{R}^{15 \times 3},$$

$$\phi_1 = M_{4,5} = [0, -e_3, -e_2] \quad \text{as well as} \quad \psi_1 = S_2 \cdot S_1^{-1}\mathcal{P}_1 + \cdots + S_5 \cdot S_4^{-1}\mathcal{P}_4 = [0, e_3, 0].$$

**Stages 6 → ∞.** Continuing the recursion, semi-infinite matrices $E$ and $M$ arise with $\phi_i = M_{4,4+i}$ and $\psi_i = S_{i+1} \cdot S_1^{-1}\mathcal{P}_1 + \cdots + S_{i+4} \cdot S_4^{-1}\mathcal{P}_4$, where in case of the particular example (3.3), it is not difficult to see by induction a periodicity of 4 within the coefficients $\phi_i$ according to

$$\phi_{1+4i} = [0, -e_3, -e_2], \quad \phi_{2+4i} = [0, e_2, 0], \quad \phi_{3+4i} = [0, 0, -(e_1 + e_2)], \quad \phi_{4+4i} = [0, e_1 + e_2, e_3]$$

for $i \geq 0$ and we obtain by use of $(1 - \varepsilon^l)^{-1} = 1 + \varepsilon^l + \varepsilon^{2l} + \cdots$ the identity

$$\phi(\varepsilon) = \begin{pmatrix} 1 & \varepsilon^4 \cdot (1 - \varepsilon^4)^{-1} & -\varepsilon^3 \cdot (1 - \varepsilon^4)^{-1} \\ 0 & (1 - \varepsilon^2)^{-1} & -\varepsilon \cdot (1 - \varepsilon^2)^{-1} \\ 0 & -\varepsilon \cdot (1 - \varepsilon^4)^{-1} & (1 - \varepsilon^4)^{-1} \end{pmatrix}.$$

Further, due to $R_4^c = \{0\}$ in (3.4), we have $R_i^c = \{0\}$ and $S_{i+1} = [0 \; 0 \; 0]$ for $i \geq 4$ and the transformation $\psi(\varepsilon)$ simplifies to a polynomial of degree 3 by

$$\psi(\varepsilon) = \begin{pmatrix} 1 & 0 & 0 \\ 0 & 1 & 0 \\ \varepsilon^3 & \varepsilon & 1 \end{pmatrix} \quad \text{with} \quad \psi^{-1}(\varepsilon) = \begin{pmatrix} 1 & 0 & 0 \\ 0 & 1 & 0 \\ -\varepsilon^3 & -\varepsilon & 1 \end{pmatrix}.$$

Hence, the families $\phi(\varepsilon)$ and $\psi^{-1}(\varepsilon)$ of the diagonalization $\Delta(\varepsilon) = \psi^{-1}(\varepsilon) \cdot L(\varepsilon) \cdot \phi(\varepsilon)$ are derived and the resolvent $L^{-1}(\varepsilon), \varepsilon \neq 0$ is conveniently calculated by the following pole of order $k = 3$

$$L^{-1}(\varepsilon) = \phi(\varepsilon) \cdot \Delta^{-1}(\varepsilon) \cdot \psi^{-1}(\varepsilon) = \begin{pmatrix} (1 - \varepsilon^4)^{-1} & 0 & -\varepsilon \cdot (1 - \varepsilon^4)^{-1} \\ (1 - \varepsilon^2)^{-1} & \varepsilon^{-2} & -\varepsilon^{-3} \cdot (1 - \varepsilon^2)^{-1} \\ -\varepsilon \cdot (1 - \varepsilon^4)^{-1} & 0 & \varepsilon^{-2} \cdot (1 - \varepsilon^4)^{-1} \end{pmatrix}.$$

The example is characterized by certain simplifications. First, $L(\varepsilon)$ is working in finite dimensions, guaranteeing stabilization of Jordain chains at some $k \geq 0$, as well closedness of all subspaces. Secondly, the matrix function $L(\varepsilon)$ is mapping between spaces of equal dimensions satisfying $N_4 = \{0\}$ and $R_4^c = \{0\}$. Thus, at $\varepsilon = 0$ an isolated singularity occurs with classical resolvent existing for $\varepsilon \neq 0$. Thirdly, $L(\varepsilon)$ is a matrix polynomial.

### 3.1 The Recursion

The proposed procedure is inherently inductive. The construction of the various operators and associated subspaces depends directly on the current state of the system. We begin by stating the assumptions and describing the stage by stage construction. When the basic procedure has been established we will use an inductive argument to justify the entire process and establish the main results.

We are starting with a formal power series of linear operators $L(\varepsilon) = \sum_{i=0}^{\infty} \varepsilon^i \cdot L_i$, $L_i \in L(B, \bar{B})$ with $B, \bar{B}$ real or complex vector spaces, $\varepsilon \in \mathbb{K} = \mathbb{R}, \mathbb{C}$ and $L(B, \bar{B})$ denoting the vector space of linear mappings between the vector spaces $B$ and $\bar{B}$. The recursion is working within this general vector space setting, where in section 5, additional assumptions are formulated to ensure local convergence of all power series involved. In this sense, we clearly distinguish what is actually needed for the recursion to work and what is needed for convergence.

The basic input to the recursion is given by the coefficients $\{L_i\}_{i \in \mathbb{N}-1}$ of $L(\varepsilon)$ and the trivial splittings $B = \{0\} \oplus B =: N_0^c \oplus N_0$ and $\bar{B} = \{0\} \oplus \bar{B} =: R_0 \oplus R_0^c$.

**Stage 1.** Define $\bar{S}_1 := L_0 \in L(B, \bar{B})$ and

$$\boxed{S_1 := I_{\bar{B}} \cdot \bar{S}_1 \in L(B, R_0^c)} \qquad \begin{aligned} N_1 &:= N[S_{1|N_0}] & N_0 &= N_1^c \oplus N_1 \\ R_1 &:= R[S_{1|N_0}] & R_0^c &= R_1 \oplus R_1^c \end{aligned} \tag{3.5}$$

with $N_1^c \subset B$ and $R_1^c \subset \bar{B}$ representing algebraic direct complements of $N_1$ in $B$ and $R_1$ in $\bar{B}$ respectively. Concerning the existence of $N_1^c$ and $R_1^c$ see [20, 4.8-A]. Thus, we obtain a split of the linear spaces $B$ and $\bar{B}$ with bijectivity of the subspaces $N_1^c$ and $R_1$ given by $S_1$. Let $\mathcal{P}_1 \in L(\bar{B}, R_1)$ denote the unique projection to $R_1$ defined by the algebraic decomposition $\bar{B} = R_1 \oplus R_1^c$ and set $E_1, M_1 \in L(B, B)$ by $E_1 = (E_{1,1}) := I_B$ and $M_1 = (M_{1,1}) := I_B$.

**Stage 2.** Define $\bar{S}_2 := L_1 \cdot M_1 \in L(B, \bar{B})$ and

$$\boxed{S_2 := (I_{\bar{B}} - \mathcal{P}_1) \cdot \bar{S}_2 \in L(B, R_1^c)} \qquad \begin{aligned} N_2 &:= N[S_{2|N_1}] & N_1 &= N_2^c \oplus N_2 \\ R_2 &:= R[S_{2|N_1}] & R_1^c &= R_2 \oplus R_2^c \end{aligned}$$

with $N_2^c$, $R_2^c$ representing algebraic direct complements of $N_2$ in $N_1$ and of $R_2$ in $R_1^c$ respectively. Let $\mathcal{P}_2 \in L(\bar{B}, R_2)$ denote the unique projection to $R_2$ defined by the algebraic decomposition $\bar{B} = R_1 \oplus R_2 \oplus R_2^c$ and set $E_2 \in L(B, B^2)$ by

$$E_2 = \begin{pmatrix} E_{1,2} \\ E_{2,2} \end{pmatrix} := \begin{pmatrix} I_B & S_1^{-1} \mathcal{P}_1 \cdot \bar{S}_2 \\ 0 & I_B \end{pmatrix}^{-1} \cdot \begin{pmatrix} 0 \\ I_B \end{pmatrix} = \begin{pmatrix} -S_1^{-1} \mathcal{P}_1 \cdot \bar{S}_2 \\ I_B \end{pmatrix}$$

with $S_1^{-1} \in L(R_1, N_1^c)$ representing the inverse of $S_1$ with respect to the subspaces $R_1$ and $N_1^c$. Finally, define $M_2 \in L(B, B^2)$ along

$$M_2 = \begin{pmatrix} M_{1,2} \\ M_{2,2} \end{pmatrix} := \begin{pmatrix} I_B & 0 \\ 0 & M_{1,1} \end{pmatrix} \cdot \begin{pmatrix} E_{1,2} \\ E_{2,2} \end{pmatrix}.$$

**Stage 3.** Define $\bar{S}_3 := [L_1, L_2] \cdot M_2 = L_1 M_{1,2} + L_2 M_{2,2} \in L(B, \bar{B})$ and

$$\boxed{S_3 := (I_{\bar{B}} - \mathcal{P}_1 - \mathcal{P}_2) \cdot \bar{S}_3 \in L(B, R_2^c)} \qquad \begin{aligned} N_3 &:= N[S_{3|N_2}] & N_2 &= N_3^c \oplus N_3 \\ R_3 &:= R[S_{3|N_2}] & R_2^c &= R_3 \oplus R_3^c \end{aligned}$$

with $N_3^c$, $R_3^c$ representing algebraic direct complements of $N_3$ in $N_2$ and of $R_3$ in $R_2^c$ respectively. Let $\mathcal{P}_3 \in L(\bar{B}, R_3)$ denote the unique projection to $R_3$ defined by the algebraic decomposition $\bar{B} = R_1 \oplus R_2 \oplus R_3 \oplus R_3^c$ and set $E_3 \in L(B, B^3)$ by

$$E_3 = \begin{pmatrix} E_{1,3} \\ E_{2,3} \\ E_{3,3} \end{pmatrix} := \begin{pmatrix} I_B & S_1^{-1}\mathcal{P}_1 \cdot \bar{S}_2 & S_1^{-1}\mathcal{P}_1 \cdot \bar{S}_3 \\ 0 & I_B & S_2^{-1}\mathcal{P}_2 \cdot \bar{S}_3 \\ 0 & 0 & I_B \end{pmatrix}^{-1} \cdot \begin{pmatrix} 0 \\ 0 \\ I_B \end{pmatrix}$$

$$= \begin{pmatrix} -S_1^{-1}\mathcal{P}_1 \cdot [I_B - \bar{S}_2 \cdot S_2^{-1}\mathcal{P}_2] \cdot \bar{S}_3 \\ -S_2^{-1}\mathcal{P}_2 \cdot \bar{S}_3 \\ I_B \end{pmatrix}$$

with $S_2^{-1} \in L(R_2, N_2^c)$ representing the inverse of $S_2$ with respect to the subspaces $R_2$ and $N_2^c$. Finally, define $M_3 \in L(B, B^3)$ along

$$M_3 = \begin{pmatrix} M_{1,3} \\ M_{2,3} \\ M_{3,3} \end{pmatrix} := \begin{pmatrix} I_B & 0 & 0 \\ 0 & M_{1,1} & M_{1,2} \\ 0 & 0 & M_{2,2} \end{pmatrix} \cdot \begin{pmatrix} E_{1,3} \\ E_{2,3} \\ E_{3,3} \end{pmatrix}.$$

The general step is described as follows.

**Stage $k + 1$.** Note that for $k = 1$ and $k = 2$, the stages 2 and Stage 3 are repeated. Define

$$\bar{S}_{k+1} := [L_1, \ldots, L_k] \cdot M_k = L_1 M_{1,k} + \cdots + L_k M_{k,k} \in L(B, \bar{B})$$

(3.6)

$$\boxed{S_{k+1} := (I_{\bar{B}} - \cdots - \mathcal{P}_k) \cdot \bar{S}_{k+1} \in L(B, R_k^c)} \qquad N_{k+1} := N[S_{k+1|N_k}] \qquad N_k = N_{k+1}^c \oplus N_{k+1}$$

$$R_{k+1} := R[S_{k+1|N_k}] \qquad R_k^c = R_{k+1} \oplus R_{k+1}^c$$

with $N_{k+1}^c, R_{k+1}^c$ representing algebraic direct complements of $N_{k+1}$ in $N_k$ and of $R_{k+1}$ in $R_{k+1}^c$ respectively. In particular, the algebraic decompositions (1.3) of $B$ and $\bar{B}$ are constructed. Let $\mathcal{P}_{k+1} \in L(\bar{B}, R_{k+1})$ denote the unique projection to $R_{k+1}$ defined by $\bar{B} = R_1 \oplus \cdots \oplus R_{k+1} \oplus R_{k+1}^c$ and set $E_{k+1} \in L(B, B^{k+1})$ by

$$E_{k+1} = \begin{pmatrix} E_{1,k+1} \\ E_{2,k+1} \\ \vdots \\ E_{k+1,k+1} \end{pmatrix} := \begin{pmatrix} I_B & S_1^{-1}\mathcal{P}_1 \cdot \bar{S}_2 & \cdots & S_1^{-1}\mathcal{P}_1 \cdot \bar{S}_{k+1} \\ 0 & I_B & \cdots & S_2^{-1}\mathcal{P}_2 \cdot \bar{S}_{k+1} \\ \vdots & \vdots & \ddots & \vdots \\ 0 & 0 & 0 & I_B \end{pmatrix}^{-1} \cdot \begin{pmatrix} 0 \\ 0 \\ \vdots \\ I_B \end{pmatrix}$$

with components $E_{i,k+1}, i = k+1, \ldots, 1$ explicitly calculated from bottom to top according to

$$E_{k+1,k+1} := I_B \in L(B, B)$$

(3.7)

$$E_{i,k+1} := -S_i^{-1}\mathcal{P}_i \cdot \sum_{\nu=i+1}^{k+1} \bar{S}_\nu \cdot E_{\nu,k+1} \in L(B, B), \quad i = k, \ldots, 1.$$

and $S_i^{-1} \in L(R_i, N_i^c)$ representing the inverse of $S_i$ with respect to the subspaces $R_i$ and $N_i^c$. Finally, define $M_{k+1} \in L(B, B^{k+1})$ along

$$M_{k+1} = \begin{pmatrix} M_{1,k+1} \\ M_{1,k+1} \\ \vdots \\ M_{k+1,k+1} \end{pmatrix} := \begin{pmatrix} I_B & 0 & \cdots & 0 \\ 0 & M_{1,1} & \cdots & M_{1,k} \\ \vdots & \vdots & \ddots & \vdots \\ 0 & 0 & 0 & M_{k,k} \end{pmatrix} \cdot \begin{pmatrix} E_{1,k+1} \\ E_{2,k+1} \\ \vdots \\ E_{k+1,k+1} \end{pmatrix}.$$

(3.8)

In general, we define invertible upper block triangular operator matrices for $k \geq 0$

$$E^{(k+1)} := \begin{pmatrix} E_{1,1} & E_{1,2} & \cdots & E_{1,k+1} \\ 0 & E_{2,2} & \cdots & E_{2,k+1} \\ \vdots & \vdots & \ddots & \vdots \\ 0 & 0 & 0 & E_{k+1,k+1} \end{pmatrix} \text{ and } M^{(k+1)} := \begin{pmatrix} M_{1,1} & M_{1,2} & \cdots & M_{1,k+1} \\ 0 & M_{2,2} & \cdots & M_{2,k+1} \\ \vdots & \vdots & \ddots & \vdots \\ 0 & 0 & 0 & M_{k+1,k+1} \end{pmatrix} \quad (3.9)$$

and we use the abbreviations $E^{(k+1)} = trng[E_1 \ldots E_{k+1}]$ and $M^{(k+1)} = trng[M_1 \ldots M_{k+1}]$ for a quadratic upper triangular matrix defined by operator column vectors of increasing length [12].

We stress that the recursion can be continued arbitrarily to define semi-infinite operator matrices $E := E^{(\infty)}$ and $M := M^{(\infty)}$. However, if the recursion stabilizes at $k \geq 0$, then the recursion simplifies significantly, due to $N_{i+1}^c = \{0\}$, $R_{i+1} = \{0\}$ for $i \geq k + 1$. As a consequence, starting with row $k + 1$, the matrix $M$ turns into a Toeplitz matrix, as will be seen below, and it is exactly row $k + 1$ of the matrix $M$ that defines the power series $\phi(\varepsilon) = \sum_{i=0}^{\infty} \varepsilon^i \cdot \phi_i$ by $\phi_i = M_{k+1,k+1+i}$, whereas the power series $\psi(\varepsilon) = \sum_{i=0}^{\infty} \varepsilon^i \cdot \psi_i$ is given by $\psi_i = S_{i+1} \cdot S_1^{-1} \mathcal{P}_1 + \cdots + S_{k+i+1} \cdot S_{k+1}^{-1} \mathcal{P}_{k+1}$. Thus, the coefficients of the transformations arise automatically from the extended recursion. Periodicity of $\phi_i$, as occurred in the example of the last section, is not a general property of the recursion.

### 3.2 Properties of Jordan chains

In the next step, let us prove some properties of the recursion concerning Jordan chains, which are summarized in Lemma 1 below. First, we easily show for $k \geq 0$ by induction

$$N[\Delta^{k+1}] = R[M^{(k+1)}{}_{|N_1 \times \cdots \times N_{k+1}}]. \quad (3.10)$$

For $k = 0$, we have $N[\Delta^1] = N[L_0] = R[I_{B|N_1}] = R[M^{(1)}{}_{|N_1}]$ and (3.10) is satisfied. Then, assume $N[\Delta^k] = R[M^{(k)}{}_{|N_1 \times \cdots \times N_k}]$ for $k \geq 0$. For calculation of $N[\Delta^{k+1}]$, we have to plug $N[\Delta^k]$ into the $k$-th leading coefficient of the Cauchy product in (3.1) and equate it to zero according to

$$(L_0 \cdots L_k) \cdot \begin{pmatrix} b_k \\ \vdots \\ b_0 \end{pmatrix} = L_0 \cdot b_k + (L_1 \cdots L_k) \cdot M^{(k)} \cdot \begin{pmatrix} n_1 \\ \vdots \\ n_k \end{pmatrix}$$

$$= L_0 \cdot n_0 + (L_1 \cdots L_k) \cdot trng[M_1, \ldots, M_k] \cdot \begin{pmatrix} n_1 \\ \vdots \\ n_k \end{pmatrix} = (\bar{S}_1 \cdots \bar{S}_{k+1}) \cdot \begin{pmatrix} n_0 \\ \vdots \\ n_k \end{pmatrix} = 0 \quad (3.11)$$

under consideration of (3.6), (3.8) and $(b_k \cdots b_0) \in B \times N[\Delta^k]$, $n_0 := b_k \in B$, $n_i \in N_i$, $i = 1, \ldots, k$. By the decomposition of $\bar{B}$ after Stage $k + 1$ (compare (1.3)), equation (3.11) is satisfied iff

$$\mathcal{P}_i \cdot (\bar{S}_1 \cdots \bar{S}_{k+1}) \cdot \begin{pmatrix} n_0 \\ \vdots \\ n_k \end{pmatrix} = 0, i = 1, \ldots, k+1 \ \wedge \ (I_{\bar{B}} - \mathcal{P}_1 - \cdots - \mathcal{P}_{k+1}) \cdot (\bar{S}_1 \cdots \bar{S}_{k+1}) \cdot \begin{pmatrix} n_0 \\ \vdots \\ n_k \end{pmatrix} = 0$$

(3.12)

$$\Leftrightarrow \begin{pmatrix} S_1 & \mathcal{P}_1 \bar{S}_2 & \cdots & \mathcal{P}_1 \bar{S}_{k+1} \\ & \ddots & \ddots & \vdots \\ & & S_k & \mathcal{P}_k \bar{S}_{k+1} \\ & & & S_{k+1} \end{pmatrix} \cdot \begin{pmatrix} n_0 \\ \vdots \\ n_{k-1} \\ n_k \end{pmatrix} = 0 \ \wedge \ \bar{S}_1 n_0 + \cdots + \bar{S}_{k+1} n_k \in R_1 \oplus \cdots \oplus R_{k+1}$$

where (3.6) implies $\bar{S}_1 n_0 \in R_1, \ldots, \bar{S}_{k+1} n_k \in R_1 \oplus \cdots \oplus R_{k+1}$ and we can restrict to the operator matrix equation in (3.12). Then, by the definitions of the operators $E_{i,j}$ in (3.7), the following equivalences result from bottom up solution of the triangular matrix system in (3.12).

$$\Leftrightarrow \begin{cases} n_k = \bar{n}_{k+1} = E_{k+1,k+1} \cdot \bar{n}_{k+1}, & \bar{n}_{k+1} \in N_{k+1} \\ n_{k-1} = \bar{n}_k - S_k^{-1} \mathcal{P}_k \bar{S}_{k+1} \cdot \bar{n}_{k+1} = \begin{bmatrix} E_{k,k} & E_{k,k+1} \end{bmatrix} \cdot \begin{pmatrix} \bar{n}_k \\ \bar{n}_{k+1} \end{pmatrix}, & \bar{n}_k \in N_k \\ \vdots & \vdots & \vdots & \vdots \end{cases}$$

$$\Leftrightarrow \begin{pmatrix} n_0 \\ \vdots \\ n_k \end{pmatrix} = \begin{pmatrix} E_{1,1} & \cdots & E_{1,k+1} \\ & \ddots & \vdots \\ & & E_{k+1,k+1} \end{pmatrix} \cdot \begin{pmatrix} \bar{n}_1 \\ \vdots \\ \bar{n}_{k+1} \end{pmatrix}, \quad \begin{pmatrix} \bar{n}_1 \\ \vdots \\ \bar{n}_{k+1} \end{pmatrix} \in N_1 \times \cdots \times N_{k+1}. \quad (3.13)$$

Thus, from (3.8), (3.11) we obtain the equivalences

$$(b_k \cdots b_0) \in N[\Delta^{k+1}] \quad (3.14)$$

$$\Leftrightarrow \begin{pmatrix} b_k \\ \vdots \\ b_0 \end{pmatrix} = \begin{pmatrix} I_B & 0 \\ 0 & M^{(k)} \end{pmatrix} \cdot \begin{pmatrix} n_0 \\ \vdots \\ n_k \end{pmatrix} = \begin{pmatrix} I_B & 0 \\ 0 & M^{(k)} \end{pmatrix} \cdot \begin{pmatrix} E_{1,1} & \cdots & E_{1,k+1} \\ & \ddots & \vdots \\ & & E_{k+1,k+1} \end{pmatrix} \cdot \begin{pmatrix} \bar{n}_1 \\ \vdots \\ \bar{n}_{k+1} \end{pmatrix}$$

$$= \begin{pmatrix} I_B & 0 & \cdots & 0 \\ & M_{1,1} & \cdots & M_{1,k} \\ & & \ddots & \vdots \\ & & & M_{k,k} \end{pmatrix} \cdot \begin{pmatrix} E_{1,1} & E_{1,2} & \cdots & E_{1,k+1} \\ & E_{2,2} & \cdots & E_{2,k+1} \\ & & \ddots & \vdots \\ & & & E_{k+1,k+1} \end{pmatrix} \cdot \begin{pmatrix} \bar{n}_1 \\ \vdots \\ \bar{n}_{k+1} \end{pmatrix}$$

$$= trng \left[ I_B \cdot E_1, \begin{pmatrix} I_B & 0 \\ 0 & M^{(1)} \end{pmatrix} \cdot E_2, \ldots, \begin{pmatrix} I_B & 0 \\ 0 & M^{(k)} \end{pmatrix} \cdot E_{k+1} \right] \cdot \begin{pmatrix} \bar{n}_1 \\ \vdots \\ \bar{n}_{k+1} \end{pmatrix}$$

$$= trng[M_1, \ldots, M_{k+1}] \cdot \begin{pmatrix} \bar{n}_1 \\ \vdots \\ \bar{n}_{k+1} \end{pmatrix} = M^{(k+1)} \cdot \begin{pmatrix} \bar{n}_1 \\ \vdots \\ \bar{n}_{k+1} \end{pmatrix}$$

implying the following results with respect to Jordan chains of increasing length.

**Lemma 1.** Given a power series of linear operators $L(\varepsilon) = \sum_{i=0}^{\infty} \varepsilon^i \cdot L_i$, $L_i \in L(B, \bar{B})$ with $B, \bar{B}$ real or complex vector spaces. Then, the iteration (3.5)-(3.8) is well defined for $k \geq 0$ with following properties.

(i) The kernel of $\Delta^{k+1} \in L(B^{k+1}, \bar{B}^{k+1})$ can be represented by the range space of $M^{(k+1)} \in L(B^{k+1}, B^{k+1})$ restricted to $N_1 \times \cdots \times N_{k+1}$, i.e. we obtain for $k \geq 0$

$$N[\Delta^{k+1}] = N\left[\begin{pmatrix} L_0 & \cdots & L_k \\ & \ddots & \vdots \\ & & L_0 \end{pmatrix}\right] = R[M^{(k+1)}_{|N_1 \times \cdots \times N_{k+1}}]$$

$$= \begin{pmatrix} I_B \\ 0 \\ 0 \\ \vdots \\ 0 \end{pmatrix} \cdot N_1 \oplus \begin{pmatrix} M_{1,2} \\ I_B \\ 0 \\ \vdots \\ 0 \end{pmatrix} \cdot N_2 \oplus \cdots \oplus \begin{pmatrix} M_{1,k+1} \\ M_{2,k+1} \\ \vdots \\ M_{k,k+1} \\ I_B \end{pmatrix} \cdot N_{k+1}.$$

(3.15)

(ii)
$$rank(b_0) \geq k \quad \Leftrightarrow \quad b_0 \in N_k \tag{3.16}$$
$$rank(b_0) = k \quad \Leftrightarrow \quad b_0 \in N_k \backslash N_{k+1} \tag{3.17}$$
$$rank(b_0) = k \quad for \quad b_0 \in N_{k+1}^c \subset N_k, \ b_0 \neq 0 \tag{3.18}$$

(iii) Stabilization of $L(\varepsilon)$ at $k \geq 0$ is equivalent to each of the following conditions

$$N_{k+1}^c \neq \{0\}, \quad N_{k+l}^c = \{0\} \quad for \quad l \geq 2 \tag{3.19}$$
$$\Leftrightarrow \quad R_{k+1} \neq \{0\}, \quad R_{k+l} = \{0\} \quad for \quad l \geq 2 \tag{3.20}$$
$$\Leftrightarrow \quad N_k \supsetneq N_{k+1} = N_{k+l} \quad for \quad l \geq 2 \tag{3.21}$$

**Proof of Lemma 1.** (i) follows from (3.14). Concerning (3.16), $rank(b_0) \geq k$ iff $b_0$ can be extended to a $k$-tuple $(b_{k-1} \cdots b_0) \in B^k$ that lies in $N[\Delta^k]$, which is possible iff $b_0 \in N_k$ by (3.15). Next, (3.17) follows from (3.16). Concerning (3.18), assume $0 \neq b_0 \in N_{k+1}^c$ with $rank(b_0) \neq k$, then $rank(b_0) \geq k+1$, yielding $b_0 \in N_{k+1}$, which contradicts $N_k = N_{k+1}^c \oplus N_{k+1}$ from (3.6).

Concerning (iii), first assume stabilization of $L(\varepsilon)$ at $k$ with $N_{k+1}^c = \{0\}$, then $N_k = \{0\} \oplus N_{k+1}$, ensuring $rank(b_0) \geq k+1$ for all $b_0 \in N_k$. Hence, there exists no element in $b_0 \in B$ with $rank(b_0) = k$, contradicting stabilization at $k$. Further, assume $N_{k+l}^c \neq \{0\}, l \geq 2$, then we obtain $0 \neq b_0 \in N_{k+l}^c$ with $rank(b_0) = k+l-1 \geq k+1$, again contradicting stabilization at $k$.

Reversely, presuppose (3.19). Then, by $N_{k+1}^c \neq \{0\}$, there exists $b_0$ with $rank(b_0) = k$ by (3.18). Next, assume the existence of $b_0$ with $rank(b_0) = \bar{k}$ and $k+1 \leq \bar{k} < \infty$. Then, by (3.17) we obtain $b_0 \in N_{\bar{k}} \backslash N_{\bar{k}+1}$ yielding $N_{\bar{k}} = N_{\bar{k}+1}^c \oplus N_{\bar{k}+1}$ from (3.6) with $N_{\bar{k}+1}^c \neq \{0\}$, contradicting $N_{k+l}^c = \{0\}$ for $l \geq 2$ and stabilization at $k$ is assured.

The equivalence between (3.19) and (3.20) is a direct consequence of bijectivity between $N_i^c$ and $R_i$ via $S_i \in L(B, \bar{B}), i \geq 1$ from (3.6). Finally, the equivalence between (3.19) and (3.21) follows again from (3.6) according to $N_k = N_{k+1}^c \oplus N_{k+1}$ and $N_{k+l} = N_{k+l+1}^c \oplus N_{k+l+1}$ for $l \geq 1$. ∎

*4. Diagonalization of Operator Power Series*

First, a pre-transformation $p_k(\varepsilon)$ of the power series $L(\varepsilon) = \sum_{i=0}^{\infty} \varepsilon^i \cdot L_i$ into the form

$$S_k(\varepsilon) \coloneqq L(\varepsilon) \cdot p_k(\varepsilon) = S_1 + \varepsilon^1 \cdot S_2 + \cdots + \varepsilon^k \cdot S_{k+1} + \sum_{i=k+1}^{\infty} \varepsilon^i \cdot Q_{i+1} \tag{4.1}$$

is derived with operators $S_i, i = 1, \ldots, k+1$ given by (3.6). We assume that the iteration is performed up to and including (3.8) with $k \geq 0$. Stabilization is not presupposed. Then, due to $S_i \cdot N_{k+1} = 0, i = 1, \ldots, k+1$, a necessary condition for $p_k(\varepsilon)$ to satisfy (4.1) reads

$$L(\varepsilon) \cdot p_k(\varepsilon) \cdot n_{k+1} = \sum_{i=k+1}^{\infty} \varepsilon^i \cdot Q_{i+1} \cdot n_{k+1} \tag{4.2}$$

for $n_{k+1} \in N_{k+1}$, i.e. $b(\varepsilon) \coloneqq p_k(\varepsilon) \cdot n_{k+1}$ has to define an approximation of $L(\varepsilon) \cdot b = 0$ of order $k + 1$. But then, it is plausible that $p_k(\varepsilon)$ might be derived from Jordan chains of length $k + 1$, as characterized by Lemma 1. In particular, by Lemma 1 (i), approximations with root elements different from zero are simply constructed from the last column $M_{k+1}$ of matrix $M^{(k+1)}$ according to

$$b(\varepsilon) = \left(\varepsilon^k, \cdots, \varepsilon^1, 1\right) \cdot \begin{pmatrix} M_{1,k+1} \\ \vdots \\ M_{k,k+1} \\ I_B \end{pmatrix} \cdot n_{k+1} = \left(I_B + \varepsilon \cdot M_{k,k+1} + \cdots + \varepsilon^k \cdot M_{1,k+1}\right) \cdot n_{k+1} \quad (4.3)$$

and setting $p_k(\varepsilon) := I_B + \varepsilon \cdot M_{k,k+1} + \cdots + \varepsilon^k \cdot M_{1,k+1}$ implies, at least, the necessary condition (4.2). Next, we obtain

$$L(\varepsilon) \cdot p_k(\varepsilon) = \left(L_0 + \varepsilon \cdot L_1 + \cdots + \varepsilon^k \cdot L_k + \sum_{i=k+1}^{\infty} \varepsilon^i \cdot L_i\right) \cdot \left(I_B + \varepsilon \cdot M_{k,k+1} + \cdots + \varepsilon^k \cdot M_{1,k+1}\right)$$

$$= L_0 I_B + \varepsilon \cdot (L_0 \; L_1) \cdot \begin{pmatrix} M_{k,k+1} \\ I_B \end{pmatrix} + \cdots + \varepsilon^k \cdot (L_0 \cdots L_k) \begin{pmatrix} M_{1,k+1} \\ \vdots \\ M_{k,k+1} \\ I_B \end{pmatrix} + \sum_{i=k+1}^{\infty} \varepsilon^i \cdot Q_{i+1}$$

$$= S_1 + \varepsilon \cdot S_2 + \cdots + \varepsilon^k \cdot S_{k+1} + \sum_{i=k+1}^{\infty} \varepsilon^i \cdot Q_{i+1} \quad (4.4)$$

with remainder coefficients $Q_{i+1} \in L(B, \bar{B}), i \geq k + 1$ and for proving (4.1), it remains to verify the last equality in (4.4) with the first summand $L_0 I_B = S_1$ agreeing by (3.5). To prove the identity of the remaining summands, we use the following lemma.

**Lemma 2.** For $k \geq 0$

$$(L_0 \cdots L_k) \cdot M^{(k+1)} = (\bar{S}_1 \cdots \bar{S}_{k+1}) \cdot E^{(k+1)} = (S_1 \cdots S_{k+1}). \quad (4.5)$$

Assuming Lemma 2, the operator $(L_0 \cdots L_k) \cdot M_{k+1}$ equals $S_{k+1}$ and the identity concerning $S_{k+1}$ in (4.4) is also shown. Next, concerning $S_k$, we first obtain from (3.8) and Lemma 2

$$(L_0 \cdots L_{k-1}) \cdot \begin{pmatrix} M_{2,k+1} \\ \vdots \\ M_{k,k+1} \\ I_B \end{pmatrix} = (L_0 \cdots L_{k-1}) \cdot M^{(k)} \cdot \begin{pmatrix} E_{2,k+1} \\ \vdots \\ E_{k,k+1} \\ I_B \end{pmatrix} = (S_1 \cdots S_k) \cdot \begin{pmatrix} E_{2,k+1} \\ \vdots \\ E_{k,k+1} \\ I_B \end{pmatrix}. \quad (4.6)$$

Further, by the definition of $E_{i,k+1}, i = k, \ldots, 2$ in (3.7), we see $R[E_{i,k+1}] \subset N_i^c \subset N_{i-1}$ implying

$$S_1 \cdot E_{2,k+1} + \cdots + S_{k-1} \cdot E_{k,k+1} + S_k \cdot I_B = S_k \quad (4.7)$$

and the identity concerning $S_k$ in (4.4) is also shown. The remaining identities from $S_{k-1}$ down to $S_2$ follow in the same way. We only have to employ the equality

$$M_l^{(k+1)} = M^{(l)} \cdot E_l^{(l+1)} \cdot \ldots \cdot E_l^{(k+1)} \quad for \quad k \geq 1 \quad and \quad 1 \leq l \leq k \quad (4.8)$$

with $M_l^{(*)}$ and $E_l^{(*)}$ denoting the $(l \times l)$ triangular matrix composed of last $l$ rows and last $l$ columns of the matrices $M^{(*)}$ and $E^{(*)}$ respectively. Identity (4.8) follows by direct calculation from the definition of $M^{(k+1)}$ in (3.8). Hence (4.1) is shown and we turn to the proof of Lemma 2.

**Proof of Lemma 2.** From (3.8), (3.9) we obtain

$$(L_0 \cdots L_k) \cdot M^{(k+1)} = (L_0 \cdots L_k) \cdot trng[M_1, \ldots, M_{k+1}]$$

$$= (L_0 \cdot M_1 \mid \cdots \mid (L_0 \cdots L_k) \cdot M_{k+1}) \quad (4.9)$$

and the first component of the row vector satisfies Lemma 2 by $L_0 \cdot M_1 = L_0 \cdot I_B = \bar{S}_1 \cdot E_{1,1} = S_1$. Concerning the remaining components in (4.9), we use (3.6) and (3.8) to see for $1 \leq i \leq k$

$$(L_0 \cdots L_i) \cdot M_{i+1} = (L_0 \cdots L_i) \cdot \begin{pmatrix} I_B & 0 \\ 0 & M^{(i)} \end{pmatrix} \cdot E_{i+1} \tag{4.10}$$

$$= (L_0 \cdot I_B \mid L_1 \cdot M_1 \mid \cdots \mid (L_1 \cdots L_i) \cdot M_i) \cdot E_{i+1} = (\bar{S}_1 \mid \bar{S}_2 \mid \cdots \mid \bar{S}_{i+1}) \cdot E_{i+1}$$

yielding the first equality in Lemma 2. Next, we note that the iteratively defined mappings $E_{i,k+1}$ in (3.7) of the solution operator $E^{(k+1)}$ of the triangular system (3.12) may also be written for $k \geq 1$ in an explicit way according to

$$B_i := \bar{S}_i \, S_i^{-1} \mathcal{P}_i \in L(\bar{B}, \bar{B}), \quad i = 1, \ldots, k \tag{4.11}$$

and

$$E_{1,k+1} = -S_1^{-1} \mathcal{P}_1 \cdot (I_{\bar{B}} - B_2) \cdot \ldots \cdot (I_{\bar{B}} - B_k) \cdot I_{\bar{B}} \cdot \bar{S}_{k+1}$$

$$\vdots \tag{4.12}$$

$$E_{k-1,k+1} = -S_{k-1}^{-1} \mathcal{P}_{k-1} \cdot (I_{\bar{B}} - B_k) \cdot I_{\bar{B}} \cdot \bar{S}_{k+1}$$

$$E_{k,k+1} = -S_k^{-1} \mathcal{P}_k \cdot I_{\bar{B}} \cdot \bar{S}_{k+1}$$

$$E_{k+1,k+1} = I_B .$$

Upper components $E_{1,k+1}, \cdots, E_{k-1,k+1}$ only appear in case of $k \geq 2$. The formulas in (4.12) follow by direct inspection of the bottom up solution process of the triangular system (3.12). Note also that all of the components from $E_{1,k+1}$ down to $E_{k,k+1}$ are multiplied by $\bar{S}_{k+1}$ from the right, motivating for $k \geq 2$ the abbreviations

$$E_{i,k+1} = -S_i^{-1} \mathcal{P}_i \cdot (I_{\bar{B}} - B_{i+1}) \cdot \ldots \cdot (I_{\bar{B}} - B_k) \cdot I_{\bar{B}} \cdot \bar{S}_{k+1} =: e_{i,k+1} \cdot \bar{S}_{k+1}, \quad i = 1, \ldots, k-1$$

$$E_{k,k+1} = -S_k^{-1} \mathcal{P}_k \cdot I_{\bar{B}} \cdot \bar{S}_{k+1} =: e_{k,k+1} \cdot \bar{S}_{k+1} \tag{4.13}$$

with $e_{i,k+1} \in L(\bar{B}, B), i = 1, \ldots, k$. For completeness, we add $e_{1,2} := -S_1^{-1} \mathcal{P}_1$ and, obviously, we obtain for $k \geq 2$ the relation

$$\begin{pmatrix} e_{1,k+1} \\ \vdots \\ e_{k-1,k+1} \end{pmatrix} = \begin{pmatrix} e_{1,k} \\ \vdots \\ e_{k-1,k} \end{pmatrix} \cdot (I_{\bar{B}} - B_k) . \tag{4.14}$$

To complete the proof of Lemma 2, we need the following statement concerning $e_{i,k+1}, i = 1, \ldots, k$.

**Lemma 3.** For $k \geq 1$

$$(\bar{S}_1 \cdots \bar{S}_k) \cdot \begin{pmatrix} e_{1,k+1} \\ \vdots \\ e_{k,k+1} \end{pmatrix} = -(\mathcal{P}_1 + \cdots + \mathcal{P}_k) . \tag{4.15}$$

Using Lemma 3 and (3.6), the second equality of Lemma 2 results along the following lines for $k \geq 1$ (for $k = 0$, Lemma 2 is obviously true).

$$(\bar{S}_1 \cdots \bar{S}_{k+1}) \cdot E^{(k+1)} = (\bar{S}_1 \cdot E_1 \mid (\bar{S}_1\, \bar{S}_2) \cdot E_2 \mid \cdots \mid (\bar{S}_1 \cdots \bar{S}_{k+1}) \cdot E_{k+1})$$

$$= (\bar{S}_1 \cdot I_B \mid \bar{S}_1 \cdot e_{1,2} \cdot \bar{S}_2 + \bar{S}_2 \cdot I_B \mid \cdots \mid (\bar{S}_1 \cdots \bar{S}_k) \cdot \begin{pmatrix} e_{1,k+1} \\ \vdots \\ e_{k,k+1} \end{pmatrix} \cdot \bar{S}_{k+1} + \bar{S}_{k+1} \cdot I_B) \quad (4.16)$$

$$= (\bar{S}_1 \cdot I_B \mid (-\mathcal{P}_1 + I_B) \cdot \bar{S}_2 \mid \cdots \mid (-\mathcal{P}_1 - \cdots - \mathcal{P}_k + I_B) \cdot \bar{S}_{k+1}) = (S_1\, S_2\, \cdots\, S_{k+1})$$

This finishes the proof of Lemma 2. ∎

**Proof of Lemma 3.** The proof goes by induction with respect to $k \geq 1$. For $k = 1$, Lemma 3 is true by $\bar{S}_1 \cdot e_{1,2} = S_1 \cdot (-S_1^{-1}\mathcal{P}_1) = -\mathcal{P}_1$. Now, assume (4.15) valid with $k$ replaced by $k - 1$. Then, by (3.6), (4.11), (4.13) and (4.14)

$$(\bar{S}_1 \cdots \bar{S}_k) \cdot \begin{pmatrix} e_{1,k+1} \\ \vdots \\ e_{k,k+1} \end{pmatrix} = (\bar{S}_1 \cdots \bar{S}_{k-1}) \cdot \begin{pmatrix} e_{1,k} \\ \vdots \\ e_{k-1,k} \end{pmatrix} \cdot (I_{\bar{B}} - B_k) + \bar{S}_k \cdot (-S_k^{-1}\mathcal{P}_k)$$

$$= -(\mathcal{P}_1 + \cdots + \mathcal{P}_{k-1}) \cdot (I_{\bar{B}} - B_k) - \bar{S}_k S_k^{-1} \mathcal{P}_k \quad (4.17)$$

$$= -(\mathcal{P}_1 + \cdots + \mathcal{P}_{k-1}) \cdot I_{\bar{B}} + (\mathcal{P}_1 + \cdots + \mathcal{P}_{k-1}) \cdot \bar{S}_k S_k^{-1} \mathcal{P}_k - \bar{S}_k S_k^{-1} \mathcal{P}_k$$

$$= -(\mathcal{P}_1 + \cdots + \mathcal{P}_{k-1}) + (\mathcal{P}_1 + \cdots + \mathcal{P}_{k-1} - I_{\bar{B}}) \cdot \bar{S}_k S_k^{-1} \mathcal{P}_k$$

$$= -(\mathcal{P}_1 + \cdots + \mathcal{P}_{k-1}) - S_k S_k^{-1} \mathcal{P}_k = -(\mathcal{P}_1 + \cdots + \mathcal{P}_{k-1} + \mathcal{P}_k)$$

and the induction is finished. ∎

We collect the results concerning the pre-transformation $p_k(\varepsilon)$ and associated triangularization of $L(\varepsilon)$ in the following theorem.

**Theorem 1.** Given the iteration (3.5)-(3.8) up to $k \geq 0$ with associated direct sums of the vector spaces $B$ and $\bar{B}$.

$$\begin{array}{ccccccccccc}
 & & & & & & \overbrace{\phantom{N_i^c}}^{=\,N_i} & & & & \\
B & = & N_0^c & \oplus & N_1^c & \oplus & \cdots & \oplus & N_i^c & \oplus & \cdots & \oplus & N_{k+1}^c & \oplus & N_{k+1} \\
 & & \downarrow \boxed{S_0 := 0} & & \downarrow \boxed{S_1} & & & & \downarrow \boxed{S_i} & & & & \downarrow \boxed{S_{k+1}} & & \\
\bar{B} & = & R_0 & \oplus & R_1 & \oplus & \cdots & \oplus & R_i & \oplus & \cdots & \oplus & R_{k+1} & \oplus & R_{k+1}^c \\
 & & & & & & & & \underbrace{\phantom{\oplus R_i}}_{=\,R_i^c} & & & & & &
\end{array} \quad (4.18)$$

Then, the polynomial of degree $k$ with respect to $\varepsilon$

$$p_k(\varepsilon) = I_B + \varepsilon \cdot M_{k,k+1} + \cdots + \varepsilon^k \cdot M_{1,k+1} = (\varepsilon^k\ \cdots\ \varepsilon^1\ 1) \cdot M_{k+1} \quad (4.19)$$

transforms the power series $L(\varepsilon) = \sum_{i=0}^{\infty} \varepsilon^i \cdot L_i$, $L_i \in L(B, \bar{B})$ into

$$S_k(\varepsilon) := L(\varepsilon) \cdot p_k(\varepsilon) = S_1 + \varepsilon^1 \cdot S_2 + \cdots + \varepsilon^k \cdot S_{k+1} + \sum_{i=k+1}^{\infty} \varepsilon^i \cdot Q_{i+1} \quad (4.20)$$

with $k + 1$ leading coefficients $S_1, \ldots, S_{k+1} \in L(B, \bar{B})$ satisfying

$$S_i \cdot ( N_0^c + \cdots + N_{i-1}^c ) \subset R_{i-1}^c \quad \wedge \quad S_i \cdot N_i^c = R_i \quad \wedge \quad S_i \cdot N_i = 0 \, . \qquad (4.21)$$

Up to this point, the remainders $Q_{i+1}$ in (4.20) show no special structure, i.e. only partial triangularization of $L(\varepsilon)$ is achieved by $S_k(\varepsilon)$. Hence, in the next step, the pre-transformation $p_k(\varepsilon)$ has to be refined and eventually extended to a formal power series $\phi(\varepsilon) = I_B + \sum_{i=1}^{\infty} \varepsilon^i \cdot \phi_i$ in such a way that all remainder terms $Q_{i+1}, i \geq k+1$ are forced to map into $R_{k+1}^c$ and complete triangularization will be obtained. In some more detail, the power series $\phi(\varepsilon)$ will transform the operators $Q_{i+1}$ into the operators $S_{i+1}$ from (3.6) according to $L(\varepsilon) \cdot \phi(\varepsilon) = \sum_{i=0}^{\infty} \varepsilon^i \cdot S_{i+1}$ with

$$R[\, S_{i+1} \,] \subset R_{k+1}^c \quad \wedge \quad S_{i+1} \cdot N_{k+1} = 0 \qquad (4.22)$$

for $i \geq k+1$. Note that Theorem 1 is valid without presupposed stabilization, but for each $k$ a new transformation $p_k(\varepsilon) = \left( \varepsilon^k \cdots \varepsilon^1 \; 1 \right) \cdot M_{k+1}$ has to be build up from the previous transformation $p_{k-1}(\varepsilon) = \left( \varepsilon^{k-1} \cdots \varepsilon^1 \; 1 \right) \cdot M_k$, where the coefficients of consecutive polynomials $p_0(\varepsilon) = I_B$, $p_1(\varepsilon) = I_B + \varepsilon \cdot M_{1,2}, \ldots$ are successively defined by the columns $M_1, M_2, \ldots$ of the $M - Matrix$. In this sense, the $M - Matrix$ is the central object to be constructed.

$$M = \begin{pmatrix} I_B & M_{1,2} & \cdots & \cdots & M_{1,k+1} & M_{1,k+2} & M_{1,k+3} & M_{1,k+4} & \cdots \\ & I_B & \cdots & \cdots & M_{2,k+1} & M_{2,k+2} & M_{2,k+3} & M_{2,k+4} & \\ & \uparrow & \ddots & & \vdots & \vdots & \vdots & \vdots & \\ & \boxed{p_1(\varepsilon)} & & I_B & M_{k,k+1} & M_{k,k+2} & M_{k,k+3} & M_{k,k+4} & \\ & & & & \overbrace{I_B}^{\phi_0} & \overbrace{M_{k+1,k+2}}^{\phi_1} & \overbrace{M_{k+1,k+3}}^{\phi_2} & \overbrace{M_{k+1,k+4}}^{\phi_3} & \cdots \\ & & & & \uparrow & I_B & M_{k+2,k+3}^{\searrow} & M_{k+2,k+4}^{\searrow} & \\ & & & \boxed{p_k(\varepsilon)} & & \uparrow & I_B & M_{k+3,k+4}^{\searrow} & \\ & & & & & \boxed{p_{k+1}(\varepsilon)} & \uparrow & I_B & \cdots \\ & & & & & & \boxed{p_{k+2}(\varepsilon)} & & \ddots \end{pmatrix} \qquad (4.23)$$

Now, when going from $p_{k+1}(\varepsilon) = I_B + \varepsilon \cdot M_{k+1,k+2} + \cdots + \varepsilon^{k+1} \cdot M_{1,k+2}$ to $p_{k+2}(\varepsilon) = I_B + \varepsilon \cdot M_{k+2,k+3} + \cdots + \varepsilon^{k+2} \cdot M_{1,k+3}$, then in general the coefficients of common order in $\varepsilon$ do not agree, i.e. we have $M_{k+1,k+2} \neq M_{k+2,k+3}, \ldots$. However, as is shown below, if stabilization occurs at $k$, then the first order coefficients of $p_{k+1}(\varepsilon)$ and $p_{k+2}(\varepsilon)$ coincide by $M_{k+1,k+2} = M_{k+2,k+3}$, whereas higher order coefficients may still remain different. The fact that, under stabilization at $k$, the element $M_{k+2,k+3}$ is simply obtained by shift of the element top left is depicted in (4.23) by a small arrow according to $M_{k+2,k+3}^{\searrow}$.

Alternatively, we can say that with element $M_{k+1,k+2}$ of $p_{k+1}(\varepsilon)$, the coefficient of $\varepsilon^1$ has already reached its final configuration for all subsequent polynomials $p_{k+2}(\varepsilon), p_{k+3}(\varepsilon), \ldots$. This behaviour continues in the sense that with $M_{k+1,k+3}$, the coefficient of $\varepsilon^2$ attains its final configuration for the polynomials $p_{k+2}(\varepsilon), p_{k+3}(\varepsilon), \ldots$.

In summary, if stabilization occurs at $k \geq 0$, then row $k + 1$ of the matrix $M$ establishes the final configuration of $\varepsilon$-coefficients up to infinity and, in particular, all rows below row $k + 1$ are simply given by copies of row $k + 1$ shifted to the right and down, i.e. starting with row $k + 1$ an operator Toeplitz matrix arises. In case of $k = 0$, row $k + 1$ turns into the first row and the complete matrix $M$ represents a Toeplitz matrix.

Now, row $k + 1$ delivers the formal power series $\phi(\varepsilon)$ looked for according to

$$\phi(\varepsilon) = \sum_{i=0}^{\infty} \varepsilon^i \cdot \phi_i := I_B + \sum_{i=1}^{\infty} \varepsilon^i \cdot M_{k+1,k+1+i} \qquad (4.24)$$

and we see that $\phi(\varepsilon)$ corresponds to the Laurent series that is usually associated to a Toeplitz matrix with infinite dimensions, as introduced in [21, page 355]. Due to our assumption $L(\varepsilon) = \sum_{i=0}^{\infty} \varepsilon^i \cdot L_i$, the principal part of the Laurent series vanishes in (4.24). If $L(\varepsilon)$ would be generalized to a meromorphic operator family $L(\varepsilon) = \sum_{i=-p}^{\infty} \varepsilon^i \cdot L_i$ of order $p$, then the series in (4.24) would also contain negative exponents.

Let us now perform the transformation $L(\varepsilon) \cdot \phi(\varepsilon)$ under the assumption of the Toeplitz repetition pattern in (4.23). We obtain by Cauchy product

$$L(\varepsilon) \cdot \phi(\varepsilon) = \sum_{l=0}^{\infty} \varepsilon^l \cdot \sum_{i+j=l} L_i \phi_j = \sum_{l=0}^{\infty} \varepsilon^l \cdot (L_0 \cdots L_l) \cdot \begin{pmatrix} M_{k+1,k+1+l} \\ \vdots \\ M_{k+1,k+2} \\ I_B \end{pmatrix}$$

$$= \sum_{l=0}^{\infty} \varepsilon^l \cdot (L_0 \cdots L_l) \cdot \begin{pmatrix} M_{k+1,k+1+l} \\ \vdots \\ M_{k+l,k+1+l} \\ I_B \end{pmatrix} = \sum_{l=0}^{\infty} \varepsilon^l \cdot S_{l+1} \qquad (4.25)$$

where the last identity follows from (4.20) in Theorem 1 (replace $k$ by $k + l$) according to

$$S_{k+l}(\varepsilon) = L(\varepsilon) \cdot p_{k+l}(\varepsilon) = (L_0 + \cdots + \varepsilon^l \cdot L_l + \sum_{i=l+1}^{\infty} \varepsilon^i \cdot L_i) \qquad (4.26)$$

$$\cdot (I_B + \varepsilon \cdot M_{k+l,k+1+l} + \cdots + \varepsilon^l \cdot M_{k+1,k+1+l} + \cdots + \varepsilon^{k+l} \cdot M_{1,k+1+l})$$

$$= S_1 + \cdots + \varepsilon^l \cdot S_{l+1} + \cdots + \varepsilon^{k+l} \cdot S_{k+1+l} + \sum_{i=k+1+l}^{\infty} \varepsilon^i \cdot Q_{i+1}$$

and coefficient $S_{l+1}$ of $\varepsilon^l$ obviously satisfying by (4.26)

$$S_{l+1} = (L_0 \cdots L_l) \cdot \begin{pmatrix} M_{k+1,k+1+l} \\ \vdots \\ M_{k+l,k+1+l} \\ I_B \end{pmatrix} \qquad (4.27)$$

for $l \in \mathbb{N}$. Thus by (4.25), the coefficients of the transformed power series $L(\varepsilon) \cdot \phi(\varepsilon)$ agree with the $S$-mappings defined recursively by (3.6). Further, from (3.20), (3.21) and (4.21) we obtain for $l \geq k + 1$

$$S_{l+1} \cdot (N_0^c + \cdots + N_l^c) \subset R_l^c = R_{k+1}^c \quad \wedge \quad S_{l+1} \cdot N_{l+1} = S_{l+1} \cdot N_{k+1} = 0 \qquad (4.28)$$

implying (4.22) and complete triangularization of $L(\varepsilon)$ is achieved.

Hence, it remains to show the Toeplitz pattern in (4.23) in case of stabilization at $k$. Remember, matrix $M$ is defined column by column using (3.8), thereby invoking matrix $E$ of solution operators of the triangular system (3.12). Now, in case of stabilization at $k$, matrix $E$ adopts the form

$$E = \begin{pmatrix} I_B & E_{1,2} & \cdots & E_{1,k+1} & E_{1,k+2} & E_{1,k+3} & \cdots & E_{1,k+1+l} & \cdots \\ & \ddots & & \vdots & \vdots & \vdots & \vdots & \vdots & \\ & & I_B & E_{k,k+1} & E_{k,k+2} & E_{k,k+3} & \cdots & E_{k,k+1+l} & \cdots \\ & & & I_B & E_{k+1,k+2} & E_{k+1,k+3} & \cdots & E_{k+1,k+1+l} & \cdots \\ & & & & I_B & 0 & \cdots & 0 & \cdots \\ & & & & & I_B & \ddots & \vdots & \\ & & & & & & \ddots & 0 & \cdots \\ & & & & & & & I_B & \ddots \\ & & & & & & & & \ddots \end{pmatrix} \quad (4.29)$$

with zero operators below row $k+1$, due to (4.12) and $R_{k+l} = \{0\}$, $\mathcal{P}_{k+l} = 0$, $l \geq 2$ by (3.20). But then, from the definition of the columns of $M$ in (3.8), we obtain for last $l$ components of $M_{k+1+l}$

$$(4.30)$$

$$\begin{pmatrix} M^{\searrow}_{k+2,k+1+l} \\ \vdots \\ M^{\searrow}_{k+l,k+1+l} \\ M^{\searrow}_{k+1+l,k+1+l} \end{pmatrix} = \begin{pmatrix} I_B & M_{k+1,k+2} & \cdots & M_{k+1,k+l} \\ & \ddots & & \vdots \\ & & \ddots & M_{k+l-1,k+l} \\ & & & I_B \end{pmatrix} \cdot \begin{pmatrix} E_{k+2,k+1+l} \\ \vdots \\ E_{k+l,k+1+l} \\ E_{k+1+l,k+1+l} \end{pmatrix} = \begin{pmatrix} M_{k+1,k+l} \\ \vdots \\ M_{k+l-1,k+l} \\ I_B \end{pmatrix}$$

for $l \geq 2$. Hence, last $l$ components of $M_{k+1+l}$ in (4.23) are simply obtained by shifting last $l$ components of $M_{k+l}$ to the right and down and the repetition structure of $M$ is shown.

Concerning stabilization, triangularization and diagonalization of operator power series $L(\varepsilon)$ between real or complex vector spaces, we summarize the results in the following theorem.

**Theorem 2.** Assume stabilization of the iteration (3.5)-(3.8) at $k \geq 0$ and define a power series $\phi(\varepsilon) \in L(B,B)$ by row $k+1$ of matrix $M$ in (4.23) according to (4.24).

(i)   Then, $\phi(\varepsilon)$ transforms $L(\varepsilon) = \sum_{i=0}^{\infty} \varepsilon^i \cdot L_i$ by Cauchy product into the triangularization

$$S(\varepsilon) := L(\varepsilon) \cdot \phi(\varepsilon) = S_1 + \varepsilon^1 \cdot S_2 + \cdots + \varepsilon^k \cdot S_{k+1} + \sum_{l=2}^{\infty} \varepsilon^{k+l-1} \cdot S_{k+l} \quad (4.31)$$

with $S_1, \ldots, S_{k+1}$ satisfying (4.21) from Theorem 1 and $S_{k+l}$ satisfying

$$R[\, S_{k+l}\,] \subset R^c_{k+1} \quad \wedge \quad S_{k+l} \cdot N_{k+1} = 0 \quad (4.32)$$

for $l \geq 2$. In particular, we have $S_i \cdot N_{k+1} = 0$ for all $i \geq 1$.

(ii) The power series $S(\varepsilon) = \sum_{i=0}^{\infty} \varepsilon^i \cdot S_{i+1}$, $S_{i+1} \in L(B, \bar{B})$ from (i) can be factorized by a power series $\psi(\varepsilon)$ and a diagonal operator polynomial $\Delta(\varepsilon)$ according to

$$S(\varepsilon) = \psi(\varepsilon) \cdot \Delta(\varepsilon) \coloneqq [\, I_{\bar{B}} + \sum_{i=1}^{\infty} \varepsilon^i \cdot \psi_i \,] \cdot [\, S_1 P_1 + \varepsilon \cdot S_2 P_2 + \cdots + \varepsilon^k \cdot S_{k+1} P_{k+1} \,] \quad (4.33)$$

with $\psi_i \in L(\bar{B}, \bar{B})$ and $P_1, \ldots, P_{k+1}$ denoting projections to $N_1^c, \ldots, N_{k+1}^c$ respectively.

**Proof of Theorem 2.** (i) is shown by (4.28) and (4.30). Concerning (ii), define

$$\psi(\varepsilon) \coloneqq I_{\bar{B}} + \sum_{i=1}^{\infty} \varepsilon^i \cdot \psi_i, \quad \psi_i \coloneqq S_{i+1} \cdot S_1^{-1} \mathcal{P}_1 + \cdots + S_{k+i+1} \cdot S_{k+1}^{-1} \mathcal{P}_{k+1} \quad (4.34)$$

and rewrite

$$\psi(\varepsilon) = I_{\bar{B}} + \sum_{i=1}^{\infty} \varepsilon^i \cdot S_{i+1} \cdot S_1^{-1} \mathcal{P}_1 + \cdots + \sum_{i=1}^{\infty} \varepsilon^i \cdot S_{k+i+1} \cdot S_{k+1}^{-1} \mathcal{P}_{k+1} \quad (4.35)$$

$$= I_{\bar{B}} + \varepsilon \cdot \left( \sum_{i=1}^{\infty} \varepsilon^{i-1} \cdot S_{i+1} \right) \cdot S_1^{-1} \mathcal{P}_1 + \cdots + \varepsilon \cdot \left( \sum_{i=1}^{\infty} \varepsilon^{i-1} \cdot S_{k+i+1} \right) \cdot S_{k+1}^{-1} \mathcal{P}_{k+1}$$

$$=: I_{\bar{B}} + \varepsilon \cdot r_1^c(\varepsilon) \cdot S_1^{-1} \mathcal{P}_1 + \cdots + \varepsilon \cdot r_{k+1}^c(\varepsilon) \cdot S_{k+1}^{-1} \mathcal{P}_{k+1}$$

with power series $r_i^c(\varepsilon)$ mapping to $R_i^c$ by (4.21) and (4.32) for $i = 1, \ldots, k+1$. Now, the factorization $S(\varepsilon) = \psi(\varepsilon) \cdot \Delta(\varepsilon)$ is simply shown by Cauchy product using (4.18) according to

$$\psi(\varepsilon) \cdot \Delta(\varepsilon) = [\, I_{\bar{B}} + \varepsilon \cdot r_1^c(\varepsilon) \cdot S_1^{-1} \mathcal{P}_1 + \cdots + \varepsilon \cdot r_{k+1}^c(\varepsilon) \cdot S_{k+1}^{-1} \mathcal{P}_{k+1} \,] \cdot [S_1 P_1 + \cdots + \varepsilon^k \cdot S_{k+1} P_{k+1}]$$

$$= I_{\bar{B}} \cdot [\, S_1 P_1 + \cdots + \varepsilon^k \cdot S_{k+1} P_{k+1} \,] + \varepsilon \cdot r_1^c(\varepsilon) \cdot P_1 + \cdots + \varepsilon^{k+1} \cdot r_{k+1}^c(\varepsilon) \cdot P_{k+1}$$

$$= [\, S_1 + \varepsilon \cdot r_1^c(\varepsilon) \,] \cdot P_1 + \cdots + [\, \varepsilon^k \cdot S_{k+1} + \varepsilon^{k+1} \cdot r_{k+1}^c(\varepsilon) \,] \cdot P_{k+1} \quad (4.36)$$

$$= [\, S(\varepsilon) \,] \cdot P_1 + \cdots + [\, S(\varepsilon) - (S_1 + \cdots + \varepsilon^{k-1} \cdot S_k) \,] \cdot P_{k+1}$$

$$= S(\varepsilon) \cdot P_1 + \cdots + S(\varepsilon) \cdot P_{k+1} = S(\varepsilon) \cdot (P_1 + \cdots + P_{k+1}) = S(\varepsilon)$$

under consideration of $S(\varepsilon)_{|N_{k+1}} = 0$ by Theorem 2 (i). This finishes the proof of Theorem 2. ∎

### 5. Defining Equation and Diagonalization of Analytic Operator Functions

If stabilization at $k \geq 0$ is assumed, then by Theorem 2, diagonalization is possible for an operator power series $L(\varepsilon) = \sum_{i=0}^{\infty} \varepsilon^i \cdot L_i$, $L_i \in L(B, \bar{B})$ by use of the power series $\phi(\varepsilon) = I_B + \sum_{i=1}^{\infty} \varepsilon^i \cdot \phi_i$ in $B$ and $\psi(\varepsilon) = I_{\bar{B}} + \sum_{i=1}^{\infty} \varepsilon^i \cdot \psi_i$ in $\bar{B}$. In this section we show that these power series are in fact convergent, if we add analyticity of $L(\varepsilon)$ and continuity of projections associated to the direct sums (4.18) of $B$ and $\bar{B}$.

**Theorem 3.** Assume stabilization at level $k \geq 0$ of the analytic operator family $L(\varepsilon) = \sum_{i=0}^{\infty} \varepsilon^i \cdot L_i$, $L_i \in L[B, \bar{B}]$ bounded, with $B, \bar{B}$ real or complex Banach spaces and continuity of projections $P_i \in L[B, N_i^c]$ and $\mathcal{P}_i \in L[\bar{B}, R_i]$ for $i = 1, \ldots, k+1$.

(i) Then, the near identity transformation $\phi(\varepsilon)$ from (4.24) is analytic and transforms $L(\varepsilon)$ into the analytic normal form $S(\varepsilon)$ from (4.31) with properties given by (4.21) and (4.32).

(ii) Using the analytic near identity transformation $\psi(\varepsilon)$ of $\bar{B}$ from (4.34), the family $L(\varepsilon)$ is diagonalized according to $\psi^{-1}(\varepsilon) \cdot L(\varepsilon) \cdot \phi(\varepsilon) = \Delta(\varepsilon)$ with $\Delta(\varepsilon)$ from (4.33).

(iii) Within a punctured neighbourhood $\varepsilon \in U \setminus \{0\}$, kernels and ranges of $L(\varepsilon)$ are given by

$$N[\,L(\varepsilon)\,] = \phi(\varepsilon) \cdot N_{k+1} \quad and \quad R[\,L(\varepsilon)\,] = \psi(\varepsilon) \cdot [\,R_1 \oplus \cdots \oplus R_{k+1}\,]. \tag{5.1}$$

In particular, kernels $N[L(\varepsilon)], \varepsilon \neq 0$ and ranges $R[L(\varepsilon)], \varepsilon \neq 0$ are analytically embedded into the families $N(\varepsilon) = \phi(\varepsilon) \cdot N_{k+1}, \varepsilon \in U$ and $R(\varepsilon) = \psi(\varepsilon) \cdot [R_1 \oplus \cdots \oplus R_{k+1}], \varepsilon \in U$, i.e. smoothing of kernels and ranges occurs.

(iv) A smooth generalized inverse of the diagonal operator polynomial $\Delta(\varepsilon)$ reads for $\varepsilon \neq 0$

$$\Delta^{-1}(\varepsilon) = \varepsilon^{-k} \cdot S_{k+1}^{-1} \mathcal{P}_{k+1} + \cdots + S_1^{-1} \mathcal{P}_1 \tag{5.2}$$

with pole of order $k \geq 0$ at $\varepsilon = 0$. Correspondingly, the family

$$L^{-1}(\varepsilon) = \phi(\varepsilon) \cdot \Delta^{-1}(\varepsilon) \cdot \psi^{-1}(\varepsilon) \tag{5.3}$$

defines a generalized inverse of $L(\varepsilon)$, which is analytic for $\varepsilon \in U \setminus \{0\}$ with pole of order $k$ at $\varepsilon = 0$. For $\varepsilon \in U \setminus \{0\}$, the analytic families

$$L^{-1}(\varepsilon) \cdot L(\varepsilon) = \phi(\varepsilon) \cdot (P_1 + \cdots + P_{k+1}) \cdot \phi^{-1}(\varepsilon) \tag{5.4}$$

$$L(\varepsilon) \cdot L^{-1}(\varepsilon) = \psi(\varepsilon) \cdot (\mathcal{P}_1 + \cdots + \mathcal{P}_{k+1}) \cdot \psi^{-1}(\varepsilon) \,.$$

are representing projections to the subspaces $\phi(\varepsilon) \cdot [N_1^c \oplus \cdots \oplus N_{k+1}^c] \subset B$ and $R[L(\varepsilon)] = \psi(\varepsilon) \cdot [R_1 \oplus \cdots \oplus R_{k+1}] \subset \bar{B}$ respectively. The families can analytically be continued to $\varepsilon = 0$ by

$$\phi(0) \cdot (P_1 + \cdots + P_{k+1}) \cdot \phi^{-1}(0) = P_1 + \cdots + P_{k+1} \tag{5.5}$$

$$\psi(0) \cdot (\mathcal{P}_1 + \cdots + \mathcal{P}_{k+1}) \cdot \psi^{-1}(0) = \mathcal{P}_1 + \cdots + \mathcal{P}_{k+1} \,.$$

**Proof of Theorem 3.** Up to now, stabilization at $k$ implies the Toeplitz structure of $M$ in (4.23), as well as the zero operators of $E$ in (4.29). These patterns are sufficient to prove triangularization and diagonalization, as stated in Theorem 2. Yet, more structure is present in $M$ and $E$, finally allowing to prove Theorem 3. First note that in case of stabilization at $k \geq 0$, we have $R_i = \{0\}$ for $i \geq k + 2$ by (3.20), implying $I_{\bar{B}} - B_i = I_{\bar{B}} - \bar{S}_i S_i^{-1} \mathcal{P}_i = I_{\bar{B}}$ and by (4.14) we obtain

$$\begin{pmatrix} e_{1,k+1+l} \\ \vdots \\ e_{k+1,k+1+l} \end{pmatrix} = \begin{pmatrix} e_{1,k+2} \\ \vdots \\ e_{k+1,k+2} \end{pmatrix} \cdot (I_{\bar{B}} - B_{k+2}) \cdot \ldots \cdot (I_{\bar{B}} - B_{k+l}) = \begin{pmatrix} e_{1,k+2} \\ \vdots \\ e_{k+1,k+2} \end{pmatrix} \tag{5.6}$$

for $l \geq 2$. Hence, first $k + 1$ components of the columns $E_{k+2}, E_{k+3}, \ldots$ in (4.29) are almost equal in the sense that they differ only by multiplication with $\bar{S}_{k+2}, \bar{S}_{k+3}, \ldots$ according to

$$\begin{pmatrix} E_{1,k+2} \\ \vdots \\ E_{k+1,k+2} \end{pmatrix} = \begin{pmatrix} e_{1,k+2} \\ \vdots \\ e_{k+1,k+2} \end{pmatrix} \cdot \bar{S}_{k+2}, \quad \begin{pmatrix} E_{1,k+3} \\ \vdots \\ E_{k+1,k+3} \end{pmatrix} = \begin{pmatrix} e_{1,k+2} \\ \vdots \\ e_{k+1,k+2} \end{pmatrix} \cdot \bar{S}_{k+3}, \quad \cdots \quad . \tag{5.7}$$

Combining with zeros in (4.29), the following consequences result upon constructing the columns of the matrix $M$. For $k \geq 0$ and $l \geq 1$

$$M_{k+1+l} = \begin{pmatrix} I_B & 0 \\ 0 & M^{(k+l)} \end{pmatrix} \cdot E_{k+1+l} = \begin{pmatrix} I_B & 0 \\ 0 & M^{(k+l)} \end{pmatrix} \cdot \begin{pmatrix} E_{1,k+1+l} \\ \vdots \\ E_{k+1,k+1+l} \\ 0 \\ \vdots \\ 0 \\ I_B \end{pmatrix}$$

$$= \begin{pmatrix} \begin{pmatrix} I_B & 0 \\ 0 & M^{(k)} \end{pmatrix} \cdot \begin{pmatrix} E_{1,k+1+l} \\ \vdots \\ E_{k+1,k+1+l} \end{pmatrix} \\ \text{------------} \\ 0 \\ \vdots \\ 0 \end{pmatrix} + \begin{pmatrix} 0 \\ M_{k+l} \end{pmatrix} \cdot I_B \qquad (5.8)$$

$$= \begin{pmatrix} \begin{pmatrix} I_B & 0 \\ 0 & M^{(k)} \end{pmatrix} \cdot \begin{pmatrix} e_{1,k+2} \\ \vdots \\ e_{k+1,k+2} \end{pmatrix} \cdot \bar{S}_{k+1+l} \\ \text{----------------} \\ 0 \\ \vdots \\ 0 \end{pmatrix} + \begin{pmatrix} 0 \\ M_{k+l} \end{pmatrix} =: \begin{pmatrix} \bar{H} \\ \Theta_l \end{pmatrix} \cdot \bar{S}_{k+1+l} + \begin{pmatrix} 0 \\ M_{k+l} \end{pmatrix}$$

with $\bar{H} = (H_1 \cdots H_{k+1})^T$, $H_i \in L(\bar{B}, B)$ and $\bar{H}$ independent of $l \geq 1$. The abbreviation $\Theta_l$ is used for the operator zero column vector with $l$ components.

In the next step, the columns $M_{k+1+l}$ are restricted to first $k+1$ components and we concentrate on the columns to the right of $M_{2k+1}$ to obtain the following representation.

**Lemma 4.** Assume stabilization of the power series $L(\varepsilon)$ at $k \geq 0$. Then, for $l \geq k+1$

$$\begin{pmatrix} M_{1,k+1+l} \\ \vdots \\ M_{k+1,k+1+l} \end{pmatrix} = \begin{pmatrix} H_1 & & \\ \vdots & \ddots & \\ H_{k+1} & \cdots & H_1 \end{pmatrix} \cdot \begin{pmatrix} \bar{S}_{k+1+l} \\ \vdots \\ \bar{S}_{1+l} \end{pmatrix} =: H \cdot \begin{pmatrix} \bar{S}_{k+1+l} \\ \vdots \\ \bar{S}_{1+l} \end{pmatrix}. \qquad (5.9)$$

**Proof of Lemma 4.** For $l = k+1$ we obtain from (5.8) backwards step-by-step

$$M_{2k+2} \stackrel{l=k+1}{=} \begin{pmatrix} \bar{H} \\ \Theta_{k+1} \end{pmatrix} \cdot \bar{S}_{2k+2} + \begin{pmatrix} 0 \\ M_{2k+1} \end{pmatrix} \stackrel{l=k}{=} \begin{pmatrix} \bar{H} \\ \Theta_{k+1} \end{pmatrix} \cdot \bar{S}_{2k+2} + \left( \begin{pmatrix} \bar{H} \\ \Theta_k \end{pmatrix} \cdot \bar{S}_{2k+1} + \begin{pmatrix} 0 \\ M_{2k} \end{pmatrix} \right)$$

$$\stackrel{l=k-1}{=} \cdots \stackrel{l=1}{=} \begin{pmatrix} \bar{H} \\ \Theta_{k+1} \end{pmatrix} \cdot \bar{S}_{2k+2} + \begin{pmatrix} 0 \\ \bar{H} \\ \Theta_k \end{pmatrix} \cdot \bar{S}_{2k+1} + \cdots + \begin{pmatrix} \Theta_k \\ \bar{H} \\ 0 \end{pmatrix} \cdot \bar{S}_{k+2} + \begin{pmatrix} \Theta_{k+1} \\ M_{k+1} \end{pmatrix} \qquad (5.10)$$

and restricting to first $k+1$ components, we see that Lemma 4 is proved for $l = k+1$. Now, we proceed by a simple induction with respect to $l \geq k+1$, thus finishing the proof of Lemma 4. ∎

Our aim is to construct a power series equation, which is solved by the power series $\phi(\varepsilon)$. Then, assuming analyticity for $L(\varepsilon)$, we will see that the equation and its unique solution $\phi(\varepsilon)$ are also analytic. First, define

$$\bar{d}(\varepsilon) := \begin{pmatrix} M_{1,2k+2} \\ \vdots \\ M_{k+1,2k+2} \end{pmatrix} + \varepsilon \cdot \begin{pmatrix} M_{1,2k+3} \\ \vdots \\ M_{k+1,2k+3} \end{pmatrix} + \cdots = \sum_{l=k+1}^{\infty} \varepsilon^{l-k-1} \cdot \begin{pmatrix} M_{1,k+1+l} \\ \vdots \\ M_{k+1,k+1+l} \end{pmatrix} \quad (5.11)$$

and hence, the transformation $\phi(\varepsilon)$ can also be written in the form

$$\phi(\varepsilon) = \sum_{i=0}^{\infty} \varepsilon^i \cdot \phi_i = I_B + \varepsilon \cdot \phi_1 + \cdots + \varepsilon^k \cdot \phi_k + \varepsilon^{k+1} \cdot d_{k+1}(\varepsilon) \quad (5.12)$$

with $d_{k+1}(\varepsilon)$ denoting component $k+1$ of the operator vector $\bar{d}(\varepsilon)$. Now, by Lemma 4

$$\bar{d}(\varepsilon) = \sum_{l=k+1}^{\infty} \varepsilon^{l-k-1} \cdot \begin{pmatrix} H_1 & & \\ \vdots & \ddots & \\ H_{k+1} & \cdots & H_1 \end{pmatrix} \cdot \begin{pmatrix} \bar{S}_{k+1+l} \\ \vdots \\ \bar{S}_{1+l} \end{pmatrix} = H \cdot \sum_{l=k+1}^{\infty} \varepsilon^{l-k-1} \cdot \begin{pmatrix} \bar{S}_{k+1+l} \\ \vdots \\ \bar{S}_{1+l} \end{pmatrix} \quad (5.13)$$

$$= H \cdot \left[ \begin{pmatrix} \boxed{\bar{S}_{2k+2}} \\ \bar{S}_{2k+1} \\ \vdots \\ \bar{S}_{k+3} \\ \bar{S}_{k+2} \end{pmatrix} + \varepsilon \cdot \begin{pmatrix} \bar{S}_{2k+3} \\ \boxed{\bar{S}_{2k+2}} \\ \vdots \\ \bar{S}_{k+4} \\ \bar{S}_{k+3} \end{pmatrix} + \cdots + \varepsilon^{k-1} \cdot \begin{pmatrix} \bar{S}_{2k+1+k} \\ \bar{S}_{2k+k} \\ \vdots \\ \boxed{\bar{S}_{2k+2}} \\ \bar{S}_{2k+1} \end{pmatrix} + \varepsilon^k \cdot \begin{pmatrix} \bar{S}_{2k+2+k} \\ \bar{S}_{2k+1+k} \\ \vdots \\ \bar{S}_{2k+3} \\ \boxed{\bar{S}_{2k+2}} \end{pmatrix} + \cdots \right]$$

where each component within the square brackets contains the same power series $c(\varepsilon) := \bar{S}_{2k+2} + \sum_{i=1}^{\infty} \varepsilon^i \cdot \bar{S}_{2k+2+i}$ and by (3.6), (5.11) we obtain

$$c(\varepsilon) = \bar{S}_{2k+2} + \varepsilon \cdot \bar{S}_{2k+3} + \cdots = \sum_{i=0}^{\infty} \varepsilon^i \cdot [L_1 \ldots L_{2k+1+i}] \cdot M_{2k+1+i}$$

$$= [L_1 \ldots L_{2k+1}] \cdot \begin{pmatrix} M_{1,2k+1} \\ \vdots \\ M_{k+1,2k+1} \\ \phi_{k-1} \\ \vdots \\ I_B \end{pmatrix} + \varepsilon \cdot [L_1 \ldots L_{2k+2}] \cdot \begin{pmatrix} M_{1,2k+2} \\ \vdots \\ M_{k+1,2k+2} \\ \phi_k \\ \vdots \\ I_B \end{pmatrix}$$

$$+ \varepsilon^2 \cdot [L_1 \ldots L_{2k+3}] \cdot \begin{pmatrix} M_{1,2k+3} \\ \vdots \\ M_{k+1,2k+3} \\ \phi_{k+1} \\ \phi_k \\ \vdots \\ I_B \end{pmatrix} + \varepsilon^3 \cdot [L_1 \ldots L_{2k+4}] \cdot \begin{pmatrix} M_{1,2k+4} \\ \vdots \\ M_{k+1,2k+4} \\ \phi_{k+2} \\ \phi_{k+1} \\ \phi_k \\ \vdots \\ I_B \end{pmatrix} + \cdots$$

$$= [L_1 \ldots L_{k+1}] \cdot \begin{pmatrix} M_{1,2k+1} \\ \vdots \\ M_{k+1,2k+1} \end{pmatrix} + \varepsilon \cdot [L_1 \ldots L_{k+1}] \cdot \bar{d}(\varepsilon) \quad (5.14)$$

$$+ \left( [0 \; L_{k+2} \ldots L_{2k+1}] + \varepsilon \cdot [L_{k+2} \ldots L_{2k+2}] + \cdots \right) \cdot \begin{pmatrix} \phi_k \\ \vdots \\ I_B \end{pmatrix}$$

$$+ \; \varepsilon^2 \cdot \left( L_{k+2} \cdot \phi_{k+1} \; + \; \varepsilon \cdot [L_{k+2} \; L_{k+3}] \cdot \begin{pmatrix} \phi_{k+2} \\ \phi_{k+1} \end{pmatrix} \; + \; \varepsilon^2 \cdot [L_{k+2} \; L_{k+3} \; L_{k+4}] \cdot \begin{pmatrix} \phi_{k+3} \\ \phi_{k+2} \\ \phi_{k+1} \end{pmatrix} + \cdots \right).$$

Now, the last bracket in (5.14) agrees with $(L_{k+2} + \varepsilon \cdot L_{k+3} + \cdots) \cdot d_{k+1}(\varepsilon)$ and $d_{k+1}(\varepsilon) = \sum_{i=0}^{\infty} \varepsilon^i \cdot \phi_{k+1+i}$ from (5.12). Combining (5.13) and (5.14), we arrive at

$$\bar{d}(\varepsilon) = H \cdot \begin{pmatrix} 0 \\ \bar{S}_{2k+1} \\ \vdots \\ \bar{S}_{k+3} + \varepsilon \cdot \bar{S}_{k+4} + \cdots + \varepsilon^{k-2} \cdot \bar{S}_{2k+1} \\ \bar{S}_{k+2} + \varepsilon \cdot \bar{S}_{k+3} + \cdots + \varepsilon^{k-2} \cdot \bar{S}_{2k} + \varepsilon^{k-1} \cdot \bar{S}_{2k+1} \end{pmatrix} \; + \; H \cdot \begin{pmatrix} 1 \\ \varepsilon \\ \vdots \\ \varepsilon^{k-1} \\ \varepsilon^k \end{pmatrix} \cdot c(\varepsilon)$$

$$=: \bar{p}(\varepsilon) \; + \; H \cdot \begin{pmatrix} 1 \\ \varepsilon \\ \vdots \\ \varepsilon^k \end{pmatrix} \cdot [L_1 \ldots L_{k+1}] \cdot \begin{pmatrix} M_{1,2k+1} \\ \vdots \\ M_{k+1,2k+1} \end{pmatrix} \tag{5.15}$$

$$+ \; H \cdot \begin{pmatrix} 1 \\ \varepsilon \\ \vdots \\ \varepsilon^k \end{pmatrix} \cdot \left( [0 \; L_{k+2} \ldots L_{2k+1}] \; + \; \varepsilon \cdot [L_{k+2} \ldots L_{2k+2}] + \cdots \right) \cdot \begin{pmatrix} \phi_k \\ \vdots \\ I_B \end{pmatrix}$$

$$+ \; H \cdot \begin{pmatrix} 1 \\ \varepsilon \\ \vdots \\ \varepsilon^k \end{pmatrix} \cdot \left( \varepsilon \cdot [L_1 \ldots L_{k+1}] \cdot \bar{d}(\varepsilon) \; + \; \varepsilon^2 \cdot [\underbrace{0 \; \ldots \; 0}_{k \text{ times}} \; | \; L_{k+2} + \varepsilon \cdot L_{k+3} + \cdots] \cdot \begin{pmatrix} d_1(\varepsilon) \\ \vdots \\ d_{k+1}(\varepsilon) \end{pmatrix} \right)$$

$$=: \bar{q}(\varepsilon) \; + \; H \cdot \begin{pmatrix} 1 \\ \varepsilon \\ \vdots \\ \varepsilon^k \end{pmatrix} \cdot [L_1 \ldots L_k \; | \; L_{k+1} + \varepsilon \cdot L_{k+2} + \varepsilon^2 \cdot L_{k+3} + \cdots] \cdot \varepsilon \cdot \bar{d}(\varepsilon)$$

$$=: \bar{q}(\varepsilon) \; + \; Q(\varepsilon) \cdot \varepsilon \cdot \bar{d}(\varepsilon)$$

with $\bar{q}(\varepsilon)$ an operator power series vector with $k+1$ components and $Q(\varepsilon)$ a $(k+1) \times (k+1)$ operator power series matrix. Hence, the power series vector $\bar{d}(\varepsilon)$ satisfies

$$\bar{d}(\varepsilon) = \bar{q}(\varepsilon) + \varepsilon \cdot Q(\varepsilon) \cdot \bar{d}(\varepsilon) \qquad \text{or} \qquad [I_{B^{k+1}} - \varepsilon \cdot Q(\varepsilon)] \cdot \bar{d}(\varepsilon) = \bar{q}(\varepsilon). \tag{5.16}$$

Now, the equation $[I_{B^{k+1}} - \varepsilon \cdot Q(\varepsilon)] \cdot \bar{d} = \bar{q}(\varepsilon), \bar{d} \in B^{k+1}$ turns into an analytic equation, as soon as analyticity of $\bar{q}(\varepsilon)$ and $Q(\varepsilon)$ is assured, where by (5.8), (5.9) and (5.15), $\bar{q}(\varepsilon)$ and $Q(\varepsilon)$ depend of a finite combination of operators from $M^{(2k+1)}$ and the power series

$$\begin{aligned} \bar{L}_{k+1}(\varepsilon) &:= L_{k+1} + \varepsilon \cdot L_{k+2} + \varepsilon^2 \cdot L_{k+3} + \cdots \\ &\vdots \\ \bar{L}_{2k+1}(\varepsilon) &:= L_{2k+1} + \varepsilon \cdot L_{2k+2} + \varepsilon^2 \cdot L_{2k+3} + \cdots \end{aligned} \tag{5.17}$$

which turn into analyticity, as soon as analyticity of $L(\varepsilon) = \sum_{i=0}^{\infty} \varepsilon^i \cdot L_i$ is assumed. Further, the elements of $M^{(2k+1)}$ result from the iteration (3.5)-(3.8) performed up to stage $2k+1$, hence depending from a finite combination of $2k+1$ operators $L_0, \ldots, L_{2k}$ and projections to subspaces given by the direct sums of $B$ and $\bar{B}$ in (4.18). Thus, to guarantee analyticity of $\bar{q}(\varepsilon)$ and $Q(\varepsilon)$ in (5.15), we have to assume analyticity of $L(\varepsilon)$ and continuity of the projections $P_i \in L[B, N_i^c]$ and $\mathcal{P}_i \in L[\bar{B}, R_i]$ for $i = 1 \ldots, k+1$. Note also that in (3.7), the operator $S_i^{-1} \in L(R_i, N_i^c)$ is used, representing the inverse of $S_i$ with respect to the subspaces $R_i$ and $N_i^c$. Now, if the projections $P_i$ to $N_i^c$ and $\mathcal{P}_i$ to $R_i$ are continuous, then $N_i^c$ and $R_i$ are closed and $S_i^{-1} \in L[R_i, N_i^c]$ is a continuous operator in Banach spaces by bounded inverse theorem [20, 4.2-H].

Hence, under the assumptions of Theorem 3, $\bar{q}(\varepsilon)$ and $Q(\varepsilon)$ are analytic and $I - \varepsilon \cdot Q(\varepsilon)$ is an isomorphism in $B^{k+1}$ for $\varepsilon$ chosen sufficiently small, implying analyticity of the solution vector $\bar{d}(\varepsilon)$, as well as analyticity of the diffeomorphic transformation $\phi(\varepsilon)$ in (5.12). Next, by the product $S(\varepsilon) = L(\varepsilon) \cdot \phi(\varepsilon)$, the triangularization $S(\varepsilon)$ is analytic and from (4.35) we obtain analyticity of $\psi(\varepsilon)$ under consideration of $r_l^c(\varepsilon) = \sum_{i=1}^{\infty} \varepsilon^{i-1} \cdot S_{i+l}$, $l = 1, \cdots, k+1$. Hence, (i) and (ii) of Theorem 3 are shown.

Concerning (iii), note that for $\varepsilon \neq 0$ we have $N[\Delta(\varepsilon)] \equiv N_{k+1}$, implying $L(\varepsilon) \cdot b = \psi(\varepsilon) \cdot \Delta(\varepsilon) \cdot \phi^{-1}(\varepsilon) \cdot b = 0$ iff $\phi^{-1}(\varepsilon) \cdot b \in N_{k+1}$ iff $N[L(\varepsilon)] = \phi(\varepsilon) \cdot N_{k+1}$. Again for $\varepsilon \neq 0$, the identity $R[\Delta(\varepsilon)] \equiv R_1 \oplus \cdots \oplus R_{k+1}$ is true, yielding $R[L(\varepsilon)] = R[\psi(\varepsilon) \cdot \Delta(\varepsilon) \cdot \phi^{-1}(\varepsilon)] = \psi(\varepsilon) \cdot [R_1 \oplus \cdots \oplus R_{k+1}]$ and (5.1) is true.

Concerning (iv), we see by direct calculation $\Delta^{-1}(\varepsilon) \cdot \Delta(\varepsilon) = P_1 + \cdots + P_{k+1}$ and $\Delta(\varepsilon) \cdot \Delta^{-1}(\varepsilon) = \mathcal{P}_1 + \cdots + \mathcal{P}_{k+1}$, implying $\Delta(\varepsilon) \cdot \Delta^{-1}(\varepsilon) \cdot \Delta(\varepsilon) = \Delta(\varepsilon)$ and $\Delta^{-1}(\varepsilon) \cdot \Delta(\varepsilon) \cdot \Delta^{-1}(\varepsilon) = \Delta^{-1}(\varepsilon)$ and (5.2) is shown. But then, $L^{-1}(\varepsilon) \cdot L(\varepsilon) = \phi(\varepsilon) \cdot (P_1 + \cdots + P_{k+1}) \cdot \phi^{-1}(\varepsilon)$ and $L(\varepsilon) \cdot L^{-1}(\varepsilon) = \psi(\varepsilon) \cdot (\mathcal{P}_1 + \cdots + \mathcal{P}_{k+1}) \cdot \psi^{-1}(\varepsilon)$, yielding $L(\varepsilon) \cdot L^{-1}(\varepsilon) \cdot L(\varepsilon) = L(\varepsilon)$ and $L^{-1}(\varepsilon) \cdot L(\varepsilon) \cdot L^{-1}(\varepsilon) = L^{-1}(\varepsilon)$ as well as (5.3), (5.4) and (5.5). This finishes the proof of Theorem 3. ∎

**Remark 1.** Under the assumptions of Theorem 3, continuity of the projections $P_i \in L[B, N_i^c]$ and $\mathcal{P}_i \in L[\bar{B}, R_i]$ in (4.18) is equivalent to closedness of $N_i^c$ and $R_i^c$ for $i = 1 \ldots, k+1$. In other words, if it is possible to choose closed complements $N_i^c$ and $R_i^c$ in step (3.6) of the recursion for $i = 1, \ldots, k+1$, then $\phi(\varepsilon)$ is analytic (provided $L(\varepsilon)$ is analytic and stabilizes at level $k$).

**Remark 2.** If $B$ and $\bar{B}$ are Hilbert spaces, then, under the assumptions of Theorem 3, $N_{i-1}$ and $N_i$ are closed and the complement $N_i^c$ in $N_{i-1} = N_i^c \oplus N_i$ can always be chosen orthogonal by $N_i^c \coloneqq N_i^\perp \cap N_{i-1}$ with $N_i^c$ closed [20, 4.82-A]. Hence, according to Remark 1, continuity of the projections is equivalent to closedness of $R_i^c$ in $R_{i-1}^c = R_i \oplus R_i^c$ for $i = 1, \ldots, k+1$. Now, if we require $R_i$ closed, we can choose $R_i^c \coloneqq R_i^\perp \cap R_{i-1}^c$ with $R_i^c$ closed and the projections are continuous.

**Remark 3.** If $L_0 \in L[B, \bar{B}]$ is a Fredholm operator in Banach spaces, then we can be sure that the recursion for $L(\varepsilon) = L_0 + \sum_{i=1}^{\infty} \varepsilon^i \cdot L_i$ stabilizes at some value $k \geq 0$ and, under the assumptions of Theorem 3, continuity of the projections is also satisfied. Exemplarily, this constellation occurs in the matrix case with $B = \mathbb{R}^m$ and $\bar{B} = \mathbb{R}^n$. Note also that stabilization at some level $k$ is achieved, as soon as a first operator $S_i$ appears satisfying $dim(N_i) < \infty$ or $dim(R_i^c) < \infty$, i.e. in Banach spaces a semi-Fredholm operator $S_i$ establishes stabilization at some $k \geq i - 1$.

## 6. Implications

**1.** The diagonal operator polynomial $\Delta(\varepsilon)$ may further be factorized according to $\Delta(\varepsilon) = S_1 P_1 + \cdots + \varepsilon^k \cdot S_{k+1} P_{k+1} = (S_1 P_1 + \cdots + S_{k+1} P_{k+1}) \cdot (P_1 + \cdots + \varepsilon^k \cdot P_{k+1}) =: S_P \cdot P(\varepsilon)$ with an $\varepsilon$ independent operator $S_P \in L[B, \bar{B}]$. Then (1.1) reads $\psi^{-1}(\varepsilon) \cdot L(\varepsilon) \cdot \phi(\varepsilon) = S_P \cdot P(\varepsilon)$ and factorization of $L(\varepsilon)$ to the constant operator $S_P$ is possible for $\varepsilon \neq 0$ according to

$$\psi^{-1}(\varepsilon) \cdot L(\varepsilon) \cdot \phi(\varepsilon) \cdot P^{-1}(\varepsilon) = S_P \tag{6.1}$$

with $k$-pole $P^{-1}(\varepsilon) = \varepsilon^{-k} \cdot P_{k+1} + \cdots + P_1 \in L[B, N_1^c \oplus \cdots \oplus N_{k+1}^c]$. In particular, the left hand side in (6.1) can now smoothly be continued to $\varepsilon = 0$ by $S_P$. In a nonlinear context, this kind of factorization and blow up procedure is used in [14] and [18, Theorems 1,2] to derive a regular system at $\varepsilon = 0$, suitable for direct application of the implicit function theorem.

Alternatively, we obtain from $\psi^{-1}(\varepsilon) \cdot L(\varepsilon) \cdot \phi(\varepsilon) = S_P \cdot P(\varepsilon)$ the so called Smith factorization

$$L(\varepsilon) = \psi(\varepsilon) \cdot (S_1 P_1 + \cdots + S_{k+1} P_{k+1}) \cdot (P_1 + \cdots + \varepsilon^k \cdot P_{k+1}) \cdot \phi^{-1}(\varepsilon) \tag{6.2}$$

$$=: \mathcal{A}(\varepsilon) \cdot (P_1 + \cdots + \varepsilon^k \cdot P_{k+1}) \cdot \phi^{-1}(\varepsilon)$$

with $\mathcal{A}(\varepsilon)$ an analytic family from $B$ to $\bar{B}$. For $L_0$ a Fredholm operator of arbitrary index $\geq 0$, formula (6.2) can be found in [9, Definition 11.6.1, Theorem 11.6.4]. For $L_0$ a Fredholm operator of index 0, the Smith form is derived in [12, Theorem 7.8.3] and [13, Theorem 1.8.4]. Then it is also possible to represent the resolvent $L^{-1}(\varepsilon)$ in an optimal way by use of the Keldysh theorem and the Jordan chains can effectively be calculated by contour integrals as demonstrated in [7].

**2.** In [1, section 7], an operator polynomial $L(\varepsilon) = L_0 + \cdots + \varepsilon^n \cdot L_n \in L[B, \bar{B}]$ of degree $n \geq 1$ is related to its augmented linear operator pencil $\bar{L}(\varepsilon) := \bar{L}_0 + \varepsilon \cdot \bar{L}_1 \in L[B^n, \bar{B}^n]$ using

$$\bar{L}_0 := \begin{pmatrix} L_0 & & \\ \vdots & \ddots & \\ L_{n-1} & \cdots & L_0 \end{pmatrix} \quad \text{and} \quad \bar{L}_1 := \begin{pmatrix} L_n & \cdots & L_1 \\ & \ddots & \vdots \\ & & L_n \end{pmatrix}.$$

It is shown that the resolvent $R(\varepsilon)$ of the operator polynomial $L(\varepsilon)$ is well defined iff the corresponding resolvent $\bar{R}(\varepsilon)$ of the linear extended operator pencil $\bar{L}(\varepsilon)$ is well defined. In particular, properties of $\bar{R}(\varepsilon)$ can be transferred to properties of $R(\varepsilon)$ and in this sense, it is sufficient to analyze the resolvent $\bar{R}(\varepsilon)$ of the augmented linear operator pencil $\bar{L}(\varepsilon)$.

Now, the recursion (3.5)-(3.8) can be performed for both $L(\varepsilon)$ and $\bar{L}(\varepsilon)$ with corresponding operators and decompositions of the spaces $B, \bar{B}$ and $B^n, \bar{B}^n$ respectively. Then, it is not difficult to see the equivalence

$$L(\varepsilon) \text{ stabilizes at } k \geq 0 \quad \Leftrightarrow \quad \bar{L}(\varepsilon) \text{ stabilizes at } \bar{k} \geq 0$$

with $k, \bar{k}$ satisfying $(\bar{k} - 1) \cdot n < k \leq \bar{k} \cdot n$. In particular, the Jordan chains of length $\bar{k}$ of the linear pencil $\bar{L}(\varepsilon)$ and the Jordan chains of length $k$ of the polynomial $L(\varepsilon)$ are related by

$$N[\bar{\Delta}^{\bar{k}}] = \begin{pmatrix} N[\Delta^k] \\ 0 \end{pmatrix}$$

with $0 \in B^{\bar{k} \cdot n - k}$. Other relations between the recursions of $L(\varepsilon)$ and $\bar{L}(\varepsilon)$ follow along similar lines of reasoning.

**3.** In [1], [2], linear operator pencils of the form $L(\varepsilon) = L_0 + \varepsilon \cdot L_1$, $\mathbb{K} = \mathbb{C}$, $L_0, L_1 \in L[B, \bar{B}]$, $L_0$ singular are investigated with respect to the resolvent $R(\varepsilon), \varepsilon \neq 0$. There, it is shown that if the resolvent exists as a Laurent series $R(\varepsilon) = \sum_{j \in \mathbb{Z}} \varepsilon^j \cdot R_j$, then all coefficients can be calculated recursively from the basic coefficients $\{R_{-1}, R_0\}$ according to $R_{-j} = (-1)^{j-1} \cdot (R_{-1} L_0)^{j-1} \cdot R_{-1}$ and $R_j = (-1)^j \cdot (R_0 L_1)^j \cdot R_0$ for $j \geq 1$.

In the paper at hand, we formulate sufficient conditions for $L(\varepsilon)$ at $\varepsilon = 0$, guaranteeing for $\varepsilon \neq 0$ a $k$-pole generalized inverse $L^{-1}(\varepsilon) = \phi(\varepsilon) \cdot \Delta^{-1}(\varepsilon) \cdot \psi^{-1}(\varepsilon) = \varepsilon^{-k} \cdot S_{k+1}^{-1} \mathcal{P}_{k+1} + \cdots$. In particular, $L^{-1}(\varepsilon)$ turns into a $k$-pole resolvent $R(\varepsilon)$, presupposed stabilization occurs with $N_{k+1} = \{0\}$ and $R_{k+1}^c = \{0\}$. Then, the coefficients $R_{\pm j}$ of the resolvent $R(\varepsilon)$ are expressed by operators $S_i, i \geq 1$ from (3.6) and projections to subspaces from (4.18) and it might be interesting to compare the different representations of $R_{\pm j}$.

**4.** In [14, 155-156], the recursion from section 3.1 is used to derive a sequence of bifurcation results for nonlinear operator equations of the form $G[z] = 0, G: B \to \bar{B}$ with $B, \bar{B}$ general Banach spaces. Therein, the linear operator family $L(\varepsilon)$ arises by linearization of $G$ along a curve $z(\varepsilon)$ according to $L(\varepsilon) \coloneqq G'[z(\varepsilon)] \in L[B, \bar{B}]$, where the curve $z(\varepsilon)$ represents an approximative solution curve of $G[z] = 0$ of order $2k + 1$, i.e. $G[z(\varepsilon)] = O(\varepsilon^{2k+1})$. Then, if the recursion stabilizes at $k$ with $R_{k+1}^c = \{0\}$, the implicit function theorem can be applied to ensure the existence of a solution curve $\bar{z}(\varepsilon)$ satisfying $G[\bar{z}(\varepsilon)] = 0$ near the approximation $z(\varepsilon)$. Using this approach periodic orbits emanating from degenerate Hopf points are derived in the context of singular perturbation problems [14, page 96]. The periodic solutions represent an oscillating *CAMP* signal within a cell communication system of *Dictyostelium discoideum* [14, page 5].

We note that in [14, 295-300], the recursion is not formulated with respect to the coefficients $\{L_i\}_{i \in \mathbb{N}_{-1}}$ of the linearization $L(\varepsilon) = G'[z(\varepsilon)]$, but directly with respect to derivatives $G^{(i)}[0], i \geq 1$ of the nonlinear operator $G[z]$ at $z = 0$. This does not simplify the formulas.

Based on the formulation of the recursion in section 3.1, the analysis of nonlinear operator mappings $G: B \to \bar{B}$ is continued in [15] with respect to linearization of the nonlinear mappings and with respect to Newton's Lemma, considering the order of approximation needed to obtain solutions of $G[z] = 0$; compare [15, Theorem 1, Corollaries 1,2]. In [16, Theorem 2] power series solutions of nonlinear systems of differential-algebraic equations are analyzed using the recursion.

Finally, in [17, Corollary 4] it is shown that the Milnor number $\mu$ of a simple *ADE*-singularity $G[z]$ is easily calculated by $\mu = k_1 + \cdots + k_\tau - order(G) + 1$ with $\tau \geq 1$ denoting the number of different solution curves $z_i(\varepsilon), i = 1, \ldots, \tau$ through the singularity and corresponding stabilization of $L_i(\varepsilon) = G'[z_i(\varepsilon)]$ at $k_i \geq 0$. The formula is valid, if the homogenous polynomials of all segments in the Newton polygon factorize completely with multiplicity 1. An expected generalization of the formula with respect to higher multiplicities of the segments is also stated in [17, page 36].

**Acknowledgements:** The basic ideas of the paper go back to E. Bohl and W.-J. Beyn. It was a great pleasure to work with those people. Thanks a lot to F. Vergara for support during the final preparation of the paper. Also many thanks to the reviewer for intensive mathematical support and significant improvement of the article.